\documentclass[vecphys]{svmult}

\usepackage{graphicx}        % standard LaTeX graphics tool
\usepackage{multicol}        % used for the two-column index
\usepackage[bottom]{footmisc}% places footnotes at page bottom
\usepackage{color,amssymb,amsfonts,amsmath,stmaryrd}

\graphicspath{{Figures/}}

\newcommand{\ns}{\mathbb{N}}%{\mbox{\bbold N}}

\newcommand{\tb}{\tilde b}
\newcommand{\tc}{\tilde c}

\newcommand{\cA}{\mathcal A}  
\newcommand{\cB}{\mathcal B}   %set of bargraphs
\newcommand{\cC}{\mathcal C}
\newcommand{\cD}{\mathcal D}
\newcommand{\cF}{\mathcal F}

\newcommand{\cL}{\mathcal L}

\newcommand{\cS}{\mathcal S}	%set of staircase

\newcommand{\eps}{\epsilon}
\def\emm#1,{{\em #1}}

\newcommand{\p}{polygon}
\newcommand{\pol}{polyomino}
\newcommand{\ps}{polygons}
\newcommand{\pols}{polyominoes}
\newcommand{\bg}{bargraph}
\newcommand{\bgs}{bargraphs}
\newcommand{\cc}{column-convex}
\newcommand{\gf}{generating function}
\newcommand{\gfs}{generating functions}
\newcommand{\fps}{formal power series}

\newcommand{\figref}[1]{Fig.~\ref{#1}}
\newcommand{\secref}[1]{Section~\ref{#1}}
\begin{document}
%%%%%%%%%%%%%%%%%%%%%%%%%%%%%%%%%%%%%%%%%%%%%%%%%%%%%%%%%%%%%%%%%%%%%

\title*{Exactly solved models of polyominoes and polygons}

\author{Mireille Bousquet-M\'elou\inst{1}\and
Richard Brak\inst{2}}
\institute{CNRS, LaBRI, Universit\'e Bordeaux 1\\
 351 cours de la   Lib\'eration, 33405 Talence Cedex, France\\ 
\texttt{bousquet@labri.fr}
\and Department of Mathematics and Statistics \\
The University of Melbourne\\
Parkville 3010, Melbourne, Australia\\
\texttt{r.brak@ms.unimelb.edu.au}}
%
% Use the package "url.sty" to avoid
% problems with special characters
% used in your e-mail or web address
%
 \maketitle
 
%=================

\begin{abstract}
This chapter deals with the exact enumeration of certain classes of  
self-avoiding polygons and polyominoes on the square lattice. We  
present three general approaches that apply to many classes of  
polyominoes. The common principle to all of them is a recursive  
description of the polyominoes which then translates into a functional  
equation satisfied by the generating function. The first approach  
applies to classes of polyominoes having a linear recursive structure  
and results in a rational generating function. The second approach  
applies to classes of polyominoes having an algebraic recursive  
structure and results in an algebraic generating function. The third  
approach, commonly called the Temperley method, is based on the action  
of adding a new column  to the polyominoes. We conclude by discussing  
some open questions.

\end{abstract}
%\vskip 20mm
%
%%%%%%%%%%%%%%%%%%%%%%%%%%%%%%%%%%%%%%%%%%%%%%%%%%%%%
\section{Introduction}
%%%%%%%%%%%%%%%%%%%%%%%%%%%%%%%%%%%%%%%%%%%%%%%%%%%%%

\subsection{Subclasses of polygons and polyominoes}
%%%%%%%%%%%%%%%%%%%%%%%%%%%%%%%%%%%%%%%%%%%%%%%%%%%
This chapter deals with the \emm exact enumeration, of certain classes
of  (self-avoiding) polygons and polyominoes. 
We restrict our attention to the square lattice.
As the interior of a polygon is a polyomino, we often
consider polygons as  special polyominoes. The usual
enumeration parameters are the \emm area, (the number of cells) and
the \emm perimeter,  (the length of the border). The perimeter is always
even, and often refined into the horizontal and vertical
perimeters (number of horizontal/vertical steps in the border). Given a class
$\cC$ of polyominoes, the objective  is to
determine the following \emm complete \gf, of $\cC$:
$$
C(x,y,q)= \sum_{P\in\cC} x^{hp(P)/2} y^{vp(P)/2} q^{a(P)},
$$
where $hp(P)$, $vp(P)$ and $a(P)$ respectively denote the horizontal
perimeter, the vertical perimeter and the area of $P$. This means that
the coefficient $ c(m,n,k)$ of $x^m y^n q^k$ in the series $C(x,y,q)$ is the
number of polyominoes in the class $\cC$  having horizontal perimeter
$2m$, vertical perimeter $2n $ and area $k$. Several specializations
of $C(x,y,q)$ may be of  interest, 
such as the \emm perimeter
\gf, $C(t,t,1)$, its \emm anisotropic, version $C(x,y,1)$, or the \emm area \gf, $C(1,1,q)$. From such  exact
results, one can usually derive many of the  asymptotic properties of the
polyominoes of $\cC $: for instance the asymptotic number of polyominoes
of perimeter $n$, or the
(asymptotic) average area of these polyominoes, or
even the limiting distribution of 
this area, as $n$ tends to infinity 
 (see Chapter~11). The techniques that are used to derive 
asymptotic results from exact ones are often based on complex
analysis. A remarkable survey of these techniques is provided by
Flajolet and Sedgewick's book~\cite{flaj-sedg2}.

The study of sub-classes of polyominoes 
is natural, given the immense difficulty of
the full problem (enumerate all polygons or all polyominoes). The
objective is to develop new techniques, and to push further and further
the border between solved and unsolved models.
However, several classes have  an independent interest, other than
being an approximation of the full problem. For instance, the enumeration of
\emm partitions, (Fig.~\ref{bb:fig:zoo}(e)) is relevant in number theory and in the study of the
representations of the symmetric group. The first enumerative results
on partitions date back, at least, to Euler. A full book is devoted to
them, and is completely independent 
of the enumeration of general
polyominoes~\cite{andrews}.
Another example is provided by \emm directed, polyominoes, which are
relevant  for directed percolation, but also occur in
theoretical computer science as binary  search networks~\cite{yuba}.

All these classes will be systematically defined in
Section~\ref{sec:classes}. For the moment, let us just say that most
of them are obtained by combining
 conditions of \emm convexity, and \emm directedness.,

From the perspective of subclasses  as an approximation
to  the full problem, it is natural to ask how good this approximation
is expected to be.  The answer is quite crude:  these
approximations are terrible. For a start, all the classes that
have been counted so far are \emm exponentially small, in the class of
all polygons (or polyominoes). Hence we cannot expect their
properties to reflect faithfully those of general
polygons/polyominoes. Why would the properties of a staircase polygon
(Fig.~\ref{bb:fig:zoo}(b))  be similar to those of a general self-avoiding
polygon? Indeed, the number of staircase polygons of perimeter $2n$
grows like $2^{2n} n^{-3/2}$ (up to a multiplicative constant), while
the number of general polygons is believed to be asymptotically
$\mu^{2n} n^{-5/2}$, with $\mu=2.638\ldots$~\cite{jensen-guttmann-sap}. The average width of  a
staircase polygon is clearly linear in $n$, while the width of general
polygons is conjectured
to grow like $n^{3/4}$ (see~\cite{madras-slade}). And so on! In this context, it may be a pure
coincidence that the average area of  polygons of perimeter $2n$ is
conjectured to scale as $n^{3/2}$ (see~\cite{enting-guttmann-sap}),
\emm just as it does for staircase 
polygons,. But it is also conjectured that the limit distribution of
the  area  of $2n$-step polygons  (normalized by its average value)
coincides with 
the corresponding distribution for staircase polygons,  and for other
exactly solved classes. The universality of this distribution may not
be a coincidence (see Chapter~11 or~\cite{richard-synthese} for more
references and details).

%%%%%%%%%%%%%%%%%%%%%%%%%%%%%%%%%%%%%%%%%%%%%%%%%%%
\subsection{Three general approaches}
%%%%%%%%%%%%%%%%%%%%%%%%%%%%%%%%%%%%%%%%%%%%%%%%%%%
In this chapter, we present three robust approaches that can be applied to count
many classes $\cC$ of polyominoes. The common principle of all of them is to
translate a recursive description of the polyominoes of $\cC$ into a
functional equation satisfied by the \gf\ $C(x,y,q)$. Some readers
 may prefer seeing a translation in terms of the \emm coefficients,
of $C(x,y,q)$, namely the numbers $c(m,n,k)$. This translation is
possible, but it is usually easier to work with a functional equation
than with a recurrence relation. 
The applicability of each of these
three approaches depends on whether the polyominoes of $\cC$ have, or
don't have, a certain type of recursive structure. 

The most versatile  approach is probably the third one, as it
virtually applies to any class of polyominoes having a convexity
property. It was already used by Temperley in 1956~\cite{temperley}
and is often called, in the physics literature, \emm the Temperley
approach,. However, it often produces functional equations that are
non-trivial to solve, even when the solution finally turns out to be a
simple rational or algebraic series (these terms will be defined in
Section~\ref{sec:classes} below). 
From a combinatorics point of view, it is
important to get a better understanding of the simplicity of these
series, and this is what the first two approaches provide: the first
one applies to classes $\cC$ having a \emm linear, structure, and
gives rise to rational \gfs. The second applies to classes having an \emm
algebraic, structure, and gives rise to algebraic \gfs.

We have chosen to present these three approaches because, in our
opinion, they are the most robust ones, and we want to provide
effective tools to the reader. To our knowledge, almost all the classes that
have been solved exactly can be solved using one (or several) of these
approaches. Still, certain results have been given a beautiful
combinatorial explanation via more specific techniques. Let us 
mention two tools that are often involved in those alternative
approaches. The first tool is specific to the enumeration of polygons,
and consists in studying classes of possibly self-intersecting
polygons, and then using an inclusion-exclusion principle to eliminate
the ones with self-intersections. 
This idea appears in an old paper of P\'olya~\cite{polya}
dealing with staircase polygons,
and was further exploited to count more general
polygons~\cite{feretic-festoon1,feretic-festoon2}, including in
dimensions larger than two~\cite{guttmann-prellberg,mbm-guttmann}. 
The second tool is the use
of \emm bijections, and is of course not specific to polyomino
enumeration. The idea is to describe a one-to-one correspondence
between the objects of $\cC$ and those of another class $\cD$,
having a simpler recursive structure. In this chapter, even
though we often use  encodings of polyominoes by
words, these encodings are usually very simple and do not use the full
force of bijective methods, which is clearly at work in papers
like~\cite{mbm-viennot} or~\cite{italiens-conv-dir}.

\medskip

The structure of the chapter is simple: the three approaches we discuss
are presented, and illustrated by examples, in
Sections~\ref{sec:rat},~\ref{sec:alg} and~\ref{bb:sec:layer}
respectively. A few 
open problems which we consider worth investigating are discussed in
Section~\ref{sec:questions}.  

We conclude this introduction with  definitions of various families of
polyominoes and formal power series.

%================================================
\subsection{A visit to the zoo}\label{sec:classes}
%================================================

\begin{figure}[hbt]
\begin{center}
\includegraphics[width=12cm]{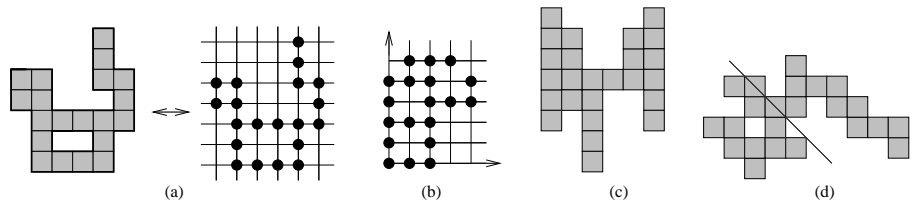}
\end{center}
\caption{From left to right: (a) a polyomino and the corresponding animal,
(b)  a NE-directed animal, (c) a column-convex polygon, (d) a $d_-$-convex polyomino.}
\label{bb:fig:poly-ani}
\end{figure}

All the classes studied in this chapter are obtained by combining
several conditions of \emm convexity, and \emm directedness., Let us
first recall that a polyomino $P$ is a finite set of square cells of the
square lattice whose interior is connected. The set of centres of the
cells form an \emm animal, $A$ (Fig.~\ref{bb:fig:poly-ani}). 
The connectivity condition means that
any two points of $A$ can be joined by a path made  up of unit
vertical and horizontal steps, in such a way that
every vertex of the path lies in $A$.
The  animal $A$  is \emm North-East directed, (or \emm directed,, for
short) if it contains a point
$v_0$, called 
the source, such that every other point of $A$ can be reached from
$v_0$ by a path made of North and East unit steps, having all its
vertices in $A$. In this case, the  polyomino corresponding to $A$ is
also said to be  NE-directed. One defines NW, SW and SE 
directed animals and polyominoes  similarly.

A polyomino $P$ is \emm column-convex,  if its
intersection with every vertical line is connected. This means that
the intersection of every vertical line with the corresponding animal
$A$ is formed 
by consecutive points. The border of $P$ is then a polygon.
 Row-convexity is defined similarly.  Finally, $P$ is
 $d_+$-convex  if the intersection of $A$ with every line of slope 1 is
 formed by consecutive points. One defines $d_-$-convex
 polyominoes  similarly. 

\begin{figure}[hbt]
\begin{center}
\includegraphics[width=12cm]{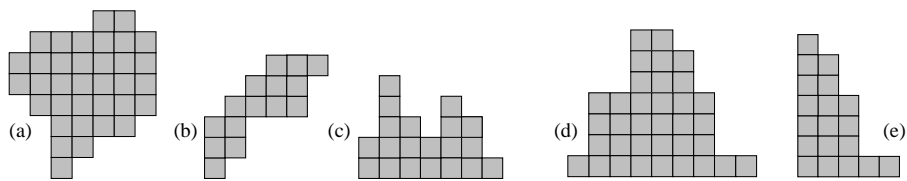}
\end{center}
\caption{A photo taken at the zoo: (a) a convex polygon, (b) a staircase
  polygon, (c) a bargraph, (d) a stack, (e) a Ferrers  diagram.}
\label{bb:fig:zoo}
\end{figure}

As discussed in~\cite{mbm-habilitation}, the combination of the 
four direction conditions and the four connectivity conditions gives rise to 31
distinct (non-symmetric) classes of polyominoes having at least one
convexity property. To these 31 classes we must add the 4 different
classes satisfying at least one directional property. 
Some prominent members of this zoo, which will occur in the forthcoming
sections,  are shown in Fig.~\ref{bb:fig:zoo}:
\begin{itemize}
\item \emm convex, polyominoes (or polygons): polyominoes that are both
  column- and row-convex,
\item \emm staircase,  polyominoes
 (or polygons): convex polygons that are NE-  and SW-directed,
\item \emm bargraphs,: column-convex polygons that are NE- and NW-directed,
\item \emm stacks,: row-convex bargraphs,
  \item \emm partitions,, a.k.a. \emm Ferrers diagrams,: convex
    polygons that are   NE-, NW- and SE-directed.
\end{itemize}

Finally, a formal power series $C(x)\equiv C(x_1, \ldots, x_k)$ with real
coefficients is \emm rational, if it
can be written as a ratio of polynomials in the $x_i$'s. It is
\emm algebraic, if it satisfies a non-trivial polynomial equation
$$
P(C(x), x_1, \ldots, x_k)=0.
$$

%%%%%%%%%%%%%%%%%%%%%%%%%%%%%%%%%%%%%%%%%%%%%%%%%%%%%%%%%%%%
\section{Linear models and rational series}
\label{sec:rat}
%%%%%%%%%%%%%%%%%%%%%%%%%%%%%%%%%%%%%%%%%%%%%%%%%%%%%%%%%%%%%

\subsection{A basic example: bargraphs counted by area}
\label{sec:bg-area}
%%%%%%%%%%%%%%%%%%%%%%%%%%%%%%%%%%%%%%%%%%%%%%%%%%%

Let $b_n$ denote the number of bargraphs of area $n$. As there is a
unique bargraph of area 1, $b_1=1$. For $n\ge 2$, there are two types
of bargraphs:
\begin{enumerate}
  \item those in which the last (\emm i.e.,, rightmost) column has height 1,
\item those in which the last  column has height 2 or more.
\end{enumerate}
Bargraphs of the first type are obtained by adding a column of height
1 to the right of any  bargraph of area $n-1$. Bargraphs of the second
type are obtained by adding one square cell to the top of the last
column of a bargraph of area $n-1$. Since a bargraph cannot be
simultaneously of type 1 and 2, this gives
$$
b_1=1 \quad \hbox{ and }  \hbox{ for } n\ge 2, \quad
b_n=2b_{n-1},
$$
which implies $b_n=2^{n-1}$. The area \gf\ of bargraphs is thus a
rational series:
$$
B(q):=\sum_{n\ge 1} b_n q^n= \frac q{1-2q}.
$$
%%%%%%%%%%%%%%%%%%%%%%%%%%%%%%%%%%%%%%%%%%%%%%%%%%%%%%%%%%%%%%%%%
\subsection{Linear objects}
%%%%%%%%%%%%%%%%%%%%%%%%%%%%%%%%%%%%%%%%%%%%%%%%%%%

The above enumeration of bargraphs is based on a very simple recursive
description 
of bargraphs. This description only involves the following two constructions:
\begin{enumerate}
  \item taking disjoint unions of sets,
\item concatenating a new cell to an already constructed object.
\end{enumerate}
In terms of \gfs\ (g.f.s),  taking the disjoint union of sets means summing their g.f.s,
while concatenating a new cell (of size 1) to all elements of a set
means multiplying its g.f.
by $q$. Hence the above description of \bgs\ translates directly into a linear
equation for the g.f.  $B(q)$:
$$
B(q)=q+qB(q)+qB(q).
$$
This equation reflects the fact that the set of bargraphs is the union
of three disjoint subsets (the unique bargraph of area 1, bargraphs of
type 1, bargraphs of type 2), and that the second and third subsets
are both obtained by adding a cell to any bargraph.

More generally, we will say that a class of objects, equipped with a
size,  is \emm linear, if
these objects can be obtained from a finite set of initial objects
using disjoint union and concatenation of one cell, or \emm
atom,. It is assumed that the concatenation of an atom increases the
size by 1. The
construction must be \emm non-ambiguous,, meaning that each object of the
class is obtained only once.  The
construction may involve  several
classes of objects simultaneously. For instance, the 
class $\tilde \cB$ of bargraphs whose last column has height 1 is linear:  the
objects of $\tilde \cB$, other than the one-cell bargraph, are
obtained by adding one cell to the right of any bargraph. The
associated series $\tilde B(q)$ is defined by the linear system:
$$
\left\{
\begin{array}{lll}
\tilde B(q)&=&q+qB(q),\\
B(q)&=&q+qB(q)+qB(q).
  \end{array}\right.
$$
In general, the \gf\ of a linear class of objects is the first
component of the solution of a system of $k$ linear equations of the form
\begin{equation}
\label{N-rational-system}
B_i(q)= P_i(q)+ q \sum_{j=1}^k a_{i,j} B_j(q)  \quad 1\le i \le k,
\end{equation}
where $a_{i,j} \in \ns$ and each $P_i(q)$ is a polynomial in $q$ with
coefficients in $\ns$. The polynomial 
$P_i(q)$ counts the initial objects of type $i$, and there are $a_{i,j}$ ways to
  aggregate an atom to an object of type $j$ to form an object of type $i$.
The system~\eqref{N-rational-system}  uniquely defines each series
$B_i(q)$, which is rational. The series obtained in this way are
called $\ns$-\emm rational,. Their study is closely related to the
theory of \emm regular languages,~\cite{salomaa}.

%%%%%%%%%%%%%%%%%%%%%%%%%%%%%%%%%%%%%%%%%%%%%%%%%%%
\subsection{More linear models}
\label{sec:linear-more}
%%%%%%%%%%%%%%%%%%%%%%%%%%%%%%%%%%%%%%%%%%%%%%%%%%%

 In this section we present three
 typical problems that can be solved
via a linear recursive description. The first one is the
perimeter enumeration of Ferrers diagrams (and stacks). The second one
generalizes the study of bargraphs performed in Section~\ref{sec:bg-area} to
all \cc\ \ps\ (and to the subclass of directed \cc\  \ps) counted by area.
The third one illustrates the role of linear models in the
approximation of  hard problems, and deals with the enumeration of
self-avoiding polygons confined to a narrow strip.
In passing, we illustrate the two following facts:  
\begin{enumerate}
  \item it may be useful to begin by describing a size preserving bijection
  between polyominoes and other objects (having a linear structure),
\item linear constructions are conveniently described by a directed
  graph when they become a bit involved.
\end{enumerate}

%%%%%%%%%%%%%%%%%%%%%%%%%%%%%%%%%%%%%%
\subsubsection{Ferrers diagrams by perimeter}
\label{bb:sec:ferrers-perimeter}
\begin{figure}[hbt]
\begin{center}
	 \scalebox{0.7}{\input{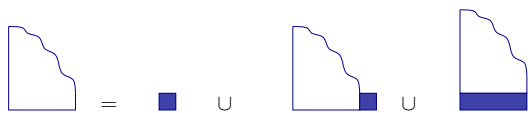}}
\end{center}
\caption{Recursive description of Ferrers diagrams.} 
\label{bb:fig:ferrers-perimeter}
\end{figure}

%%%%%%%%%%%%%%%%%%%%%%%%%%%%%%%%%%%%%%%%%%%%%%%
The set of Ferrers diagrams can be partitioned into three disjoint
subsets: first, the unique diagram of (half-)perimeter 2; then,
diagrams of width at least 2
whose rightmost column has height 1; finally, diagrams
with no column of height 1. The latter diagrams can be obtained by
duplicating the bottom row of another diagram
(Fig.~\ref{bb:fig:ferrers-perimeter}).

\begin{figure}[hbt]
\begin{center}
	 \scalebox{0.7}{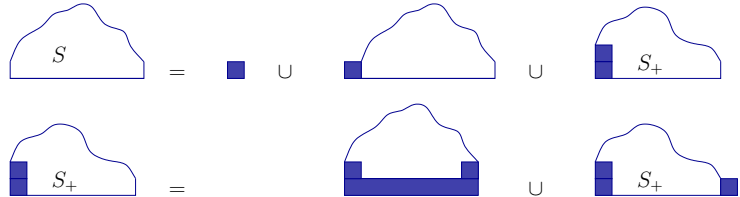}
\end{center}
\caption{Recursive description of stack polygons.} 
\label{bb:fig:stacks-perimeter}
\end{figure}

From this description, it follows that the set of words that describe
the North-East boundary of Ferrers diagrams, from the NW corner to the
SE one, admits a linear 
construction. This boundary is formed by East and South 
steps, and will be encoded by a word over the alphabet
$\{e,s\}$.  Any word over this alphabet that starts with an
$e$ and ends with an 
$s$ corresponds to a unique Ferrers diagram. Let
$\cF$ be this class of words, and let $\cL$ be the set of all non-empty prefixes of  words of $\cF$. Then $\cF$ and $\cL$
admit the following linear description:
$$
\cF=\cL s \quad \hbox{ and } \quad \cL= \{e\} \cup \cL e \cup \cL s.
$$
In these equations, the notation $\cL s$
means $\{u s, \ u\in \cL \}$, and the unions are disjoint.
 The
series that count the words of these sets by their  length (number of
letters) are thus given by the linear system
$$
F(t)=tL(t) \quad \hbox{ and } \quad  L(t)=t+2tL(t).
$$
Since the length of a coding word is the half-perimeter of the
associated diagram,
this provides the length g.f.:
$$
F(t)= \frac{t^2}{1-2t} =\sum_{n\ge 1} 2^{n-2} t^n.
$$
By separately counting  East and South steps, we obtain  the equations
\begin{equation}\label{ferrers-anisotropic}
F(x,y)=yL(x,y) \quad \hbox{and} \quad L(x,y)= x+xL(x,y)+yL(x,y),
\end{equation}
and hence the anisotropic perimeter g.f.  of these diagrams: 
$$
F(x,y)= \frac{xy}{1-x-y} =\sum_{m,n\ge 1} {{m+n-2} \choose {m-1}}x^m y^n.
$$
A similar treatment can be used to determine the perimeter g.f.  of
stack polygons: the construction  schematized in
Fig.~\ref{bb:fig:stacks-perimeter} gives:
$$ 
S(x,y)= xy+xS(x,y)+ S_+(x,y), \quad S_+(x,y)=yS(x,y)+xS_+(x,y)
$$
which yields
$$
S(x,y)= \frac{xy(1-x)}{(1-x)^2-y}.
$$

%%%%%%%%%%%%%%%%%%%%%%%%%%%%%%%%%%%%%%%%%%%%%%%%%
\subsubsection{Column-convex polygons by area}
%%%%%%%%%%%%%%%%%%%%%%%%%%%%%%%%%%%%%%%%%%%%%%%%%

\begin{figure}[hbt]
\begin{center}
\scalebox{0.7}{\input{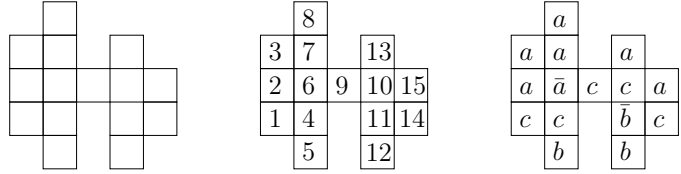}}
\end{center}
\caption{A column-convex polygon, with the numbering and encoding of
  the cells.} 
\label{bb:fig:column}
\end{figure}

Consider a \cc\ \p\ $P$  having $n$ cells. Let us number these cells from 1
to $n$ as illustrated in Fig.~\ref{bb:fig:column}. The 
columns are visited from left to right. In the first column, cells are
numbered from bottom to top. In each of the other columns, the lowest
cell that has a left neighbour gets the smallest number; then the cells
lying below it are numbered from top to bottom, and finally the cells
lying above it are numbered from bottom to top.  Note that for all $i$,
the cells labelled $1,2, \ldots , i$ form a \cc\ \p. This labelling
describes the order in which we are going to aggregate the cells.

Associate with $P$ the
word $u=u_1 \cdots u_n$ over the alphabet $\{a,b,c\}$  defined by

-- $u_i=c$ (like Column) if the $i^{\rm th}$ cell is the first to be visited in
its column,

-- $u_i=b$ (like Below) if the $i^{\rm th}$  cell lies below the first visited
cell of its column,

-- $u_i=a$ (like Above) if the $i^{\rm th}$  cell lies above the first visited
cell of its column.

\noindent Then, add a bar on the letter $u_i$ if the $i^{\rm th}$  cell of $P$
has a South
neighbour, an East neighbour, but no South-East neighbour. (In other words,
the barred letters indicate where to start a new column, when the
bottommost cell in this new column lies above the bottommost cell of
the previous column.) This gives a word $v$ over the alphabet
$\{a,b,c,\bar a , \bar b ,\bar c\}$, and $P$ can be uniquely
reconstructed from $v$.  

\begin{figure}[t]
\begin{center}
 \scalebox{0.5}{\input{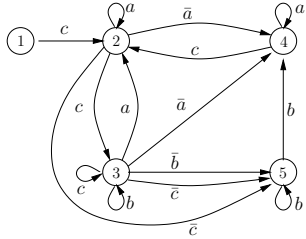}}
\end{center}
\caption{Linear construction of the words of $\cL$. The words of
  $\cL_i$ encode the paths starting from 1 and ending at $i$.}
\label{bb:fig:digraph}
\end{figure}

 We now focus on the
enumeration of these coding words.  Let $\cL$ be the
set of all prefixes
of these words, including the empty prefix $\eps$. By considering which letter can be added to the right
of which prefix, we are led to partition $\cL$ into five disjoint
subsets $\cL_1, \ldots, \cL_5$, subject to the following linear
recursive description:
\begin{equation}\label{desc:cc}
\begin{array}{llllllllllllllllllll}
{\cal L}_1& =&\{\eps\},&&&&&&&\\
{\cal L}_2& =&{\cal L}_1c \cup {\cal L}_2a \cup {\cal L}_3a\cup
{\cal L}_4c,&&&
{\cal L}_4& =&{\cal L}_2\bar a \cup {\cal L}_3\bar a \cup {\cal L}_4a \cup {\cal L}_5b,\\
{\cal L}_3& =&{\cal L}_2c\cup {\cal L}_3b\cup {\cal L}_3c,&&&
{\cal L}_5& =&{\cal L}_2 \bar c\cup{\cal L}_3\bar b\cup  {\cal L}_3\bar c
\cup {\cal L}_5b. 
\end{array}
\end{equation}
The words of ${\cal L}_4$ and ${\cal L}_5$  are those in which a barred
letter (the rightmost one) 
still waits to be ``matched'' by a letter
$c$ creating a new column. The words of $ {\cal L}_2\cup {\cal  L}_3$
are those that encode column-convex \ps.
This construction is illustrated by a directed graph in
Fig.~\ref{bb:fig:digraph}: every path starting from 1 and ending at
$i$ corresponds to a word of $\cL_i$, obtained by reading edge labels. The
series counting the words of $\cL_i$ by their length  satisfy:
$$
\begin{array}{llllllllllllllllllll}
{ L}_1& =&1,&&&&&&&\\
{ L}_2& =&q\left( { L}_1 + { L}_2 + { L}_3+{ L}_4\right),&&&
{ L}_4& =&q\left({ L}_2 + { L}_3 + { L}_4 + { L}_5\right),\\
{ L}_3& =&q\left({ L}_2+ { L}_3+ { L}_3\right),&&&
{ L}_5& =&q\left({ L}_2+{ L}_3+  { L}_3+ { L}_5\right). 
\end{array}
$$ 
The area  g.f.  of \cc\ \ps\  is $C(q)=L_2(q)+L_3(q)$. 
Solving the above system gives:
$$
C(q)={\frac {q \left( 1-q \right) ^{3}}{1-5\,q+7\,{q}^{2}-4\,{q}^{3}}}
.
$$
We believe that this result was first published by Temperley~\cite{temperley}.

A \cc \ \p\ is directed if and only if its coding word does not use the
letter $b$. We obtain a linear description of the prefixes of these
words by deleting all terms of the form $\cL_i b$ in the
description~\eqref{desc:cc}. The class $\cL_5$ becomes
irrelevant. Solving the associated system of linear equations 
gives the area g.f.  of directed \cc\ \ps:
$$
DC(q)= \frac{q(1-q)}{1-3q+q^2}.
$$
As far as we know, this result was  first published by
Klarner~\cite{klarner-results}.

%%%%%%%%%%%%%%%%%%%%%%%%%%%%%%%%%%%%%%%%%%%%%%%%%%%%%%%%%%%%%%%
\subsubsection{Polygons confined to a strip}
\label{sec:polygons-strip}
%%%%%%%%%%%%%%%%%%%%%%%%%%%%%%%%%%%%%%%%%%%%%%%%%%%%%%%%%%%%%%%

Constraining polyominoes or polygons to lie in a strip of fixed height
endows them with a linear structure. This observation 
gives a handle
to attack difficult problems, like the enumeration of general
self-avoiding polygons (SAP), self-avoiding walks, or
polyominoes~\cite{alm-janson,rote,read,zeilberger-skinny,zeilberger-animals}. As the size of the strip increases, the
approximation of the confined problem to the general one becomes better and better. This widely
applied principle gives, for instance, lower bounds on 
growth constants that are difficult to determine. We illustrate it here
with the perimeter enumeration 
of SAP confined to a strip.

Before we describe this calculation, let us mention a closely related
idea, which consists of considering anisotropic models (for
instance, SAP counted by vertical and horizontal perimeters), and 
fixing
 the number of atoms lying in one direction, for instance the
number of horizontal edges. Again, this endows the constrained objects
with a linear structure. The denominators of the rational \gfs\ that
count them often  factor in terms
$(1-y^i)$. The number of exponents $i$ that occur can be seen as a
measure of the complexity of the class. This is often observed only at an
experimental level. 
However, this observation has been pushed in some
cases to a proof that the corresponding \gf\ is not \emm D-finite,,
and in particular not algebraic (see for instance~\cite{rechni-SAP},
and Chapter~5).

\begin{figure}[hbt]
\begin{center}
	 \scalebox{0.8}{\input{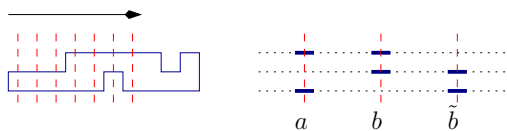}}
\end{center}
\caption{A self-avoiding polygon in a strip of height 2, encoded over a
3-letter alphabet.} 
\label{bb:fig:sap-states2}
\end{figure}

But let us return to SAP in a strip of height $k$ (a $k$\emm -strip,). A first observation
is that a polygon is completely determined by the position of its
horizontal edges.  Consider the intersection of the polygon with a
vertical line lying at a half-integer abscissa (a \emm cut,): 
the strip constraint
implies that only finitely many  configurations (or \emm states,) can occur. The
number of such states is the number of even subsets of $\{0,1, \ldots,
k\}$. This implies that SAP in a strip can be encoded by a word over a
finite alphabet. For instance, the polygon of
Fig.~\ref{bb:fig:sap-states2} is encoded by the word $\tb \tb \tb aab
aa\tb a$.

It is not hard to see that for all $k$, the set of  words encoding SAP
confined to a strip of height $k$ has a linear structure.
To make this structure clearer, we refine our encoding: for
every vertical cut, 
 we not only keep track of its
intersection with the polygon, but also of the way the horizontal
edges that  meet the cut are connected  to the left of the cut. This does
not change the size of the alphabet for $k=2$, as there is a unique
way of coupling two edges. However, if $k=3$, the configuration where
4 edges are met by the cut gives 
rise to 2 states, depending on how these 4 edges are connected
(Fig.~\ref{bb:fig:sap-states3}). The number of states is now the
number of non-crossing 
couplings on $\{0,1, \ldots, k\}$. This is also the size of our encoding
alphabet $A$.

\begin{figure}[hbt]
\begin{center}
	 \scalebox{0.7}{\input{Figures/bb-sap-states3.pstex_t}}
\end{center}
\caption{A self-avoiding polygon in a strip of height 3, encoded by
the word $d\tb a f \tc \tc e a a f \tc a c e e a b c c$.}
\label{bb:fig:sap-states3}
\end{figure}

\begin{figure}[hbt]
\begin{center}
	 \scalebox{0.5}{\input{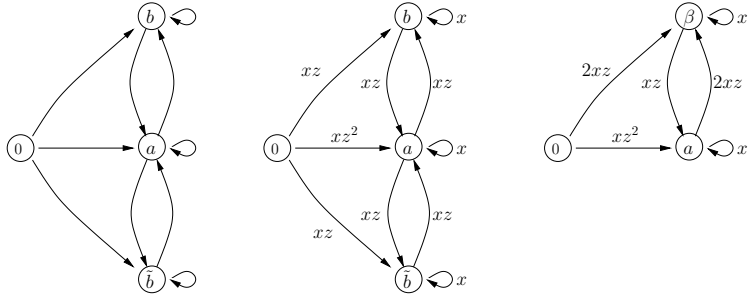}}
\end{center}
\caption{The linear structure of SAP in a 2-strip.}
\label{bb:fig:sap-automaton2}
\end{figure}

Fix $k$, and let $\cS$ be the set of words encoding SAP confined to a
$k$-strip. The 
set $\cL$ of \emm prefixes, of words of $\cS$ describes incomplete SAP, and
has a simple linear structure:  for every such prefix $w$,
the set of letters $a$ such that $wa$ lies in $\cL$ only depends on
the last letter of $w$. In other words, these prefixes are Markovian with
memory 1. For every letter $a$ in the encoding alphabet, we denote by
$\cL_a$ the set of prefixes ending with the letter $a$. The
linear structure can be encoded by a graph, from which the equations
defining the sets $\cL_a$ can automatically be written. This graph is
shown in Fig.~\ref{bb:fig:sap-automaton2} (left) for $k=2$. Every path in this graph
starting from the initial vertex $0$ corresponds to a word of
$\cL$, obtained by reading vertex labels. The linear structure of prefixes reads:
$$
\cL_a= (\epsilon + \cL_a +\cL_b +\cL_{\tb})a, \quad 
\cL_b= (\epsilon + \cL_a + \cL_b) b, \quad 
\cL_{\tb}= (\epsilon + \cL_a + \cL_{\tb}) \tb.
$$
From this we derive linear equations for incomplete SAP, where every
horizontal edge is counted by $\sqrt x$, and every vertical edge by
$\sqrt y=z$:
$$
L_a= (z^2 + L_a +zL_b +zL_{\tb})x, \quad 
L_b= (z + zL_a + L_b) x, \quad 
L_{\tb}= (z + zL_a + L_{\tb}) x.
%\end{equation}
$$
These equations keep track of how many edges are added when a new
letter is appended to a word of $\cL$. They can be schematized by a
weighted graph (Fig.~\ref{bb:fig:sap-automaton2}, middle). Now the
(multiplicative) weight 
of a path starting at $0$ is the weight of the corresponding incomplete polygon.
Finally, the completed  polygons are obtained by adding vertical edges to the
right of incomplete polygons. This gives the \gf\ of SAP in a strip of
height 2 as:
$$
S_2(x,y)= z^2L_a+ zL_b +zL_{\tb}.
$$

Clearly, we should exploit the horizontal symmetry of the model to
obtain a smaller set of equations. The letters $b$ and $\tb$ playing
symmetric roles, we replace them in the graph of
Fig.~\ref{bb:fig:sap-automaton2} by a
unique vertex $\beta$, such that the \gf\ of paths ending at $\beta$ is the
sum of the g.f.s of paths ending at $b$ and $\tb$ in the first version of the
graph (Fig.~\ref{bb:fig:sap-automaton2}, right). Introducing the series $L_\beta= L_b+L_{\tb}$,
we have thus replaced the previous system of four equations by
$$
L_a=x(z^2+L_a+zL_\beta), \quad  L_\beta=x(2z+2zL_a+L_\beta), \quad 
S_2(x,y)= z^2L_a+zL_\beta,
$$
from which we obtain
$$
S_2(x,y)=\frac{xy(2-2x+y+3xy)}{(1-x)^2-2x^2y}.
$$
Note that this series counts   polygons of height 1 twice, so that
we should subtract $S_1(x,y)=xy/(1-x)$ to obtain the g.f. of SAP of
height at most 2, defined up to translation.

\begin{figure}[hbt]
\begin{center}
	 \scalebox{0.5}{\input{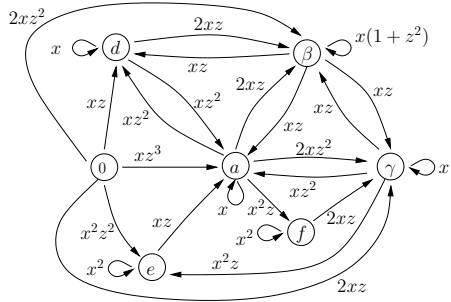}}
\end{center}
\caption{The linear structure of SAP in a 3-strip.}
\label{bb:fig:sap-automaton3}
\end{figure}

For $k=3$, the original alphabet, shown in
Fig.~\ref{bb:fig:sap-states3}, has 8 letters, 
but two pairs of them play symmetric roles. After merging the vertices
$b$ and $\tb$ on the one hand, $c$ and $\tc$ on the other, the condensed
graph, with its $x,z$ weights, is shown in
Fig.~\ref{bb:fig:sap-automaton3}. The corresponding equations read
$$
\begin{array}{l}
L_a=x \left( {z}^{3}+L_a+zL_\beta+{z}^{2}L_\gamma+{z}^{2}  L_d  +zL_e
\right) ,
%L_f={x}^{2} \left( zL_a+L_f \right) , 
\\
L_\beta=x \left( 2\,{z}^{2}+2\,zL_a+(1+z^2)L_\beta+zL_\gamma+2\,z  L_d 
 \right) ,\\
%L_e={x}^{2} \left( {z}^{2}+zL_\gamma+L_e   \right) ,\\
L_\gamma=x \left( 2\,z+2\,{z}^{2}L_a+zL_\beta+L_\gamma+2\,zL_f \right) ,\\
L_d=x \left( z+{z}^{2}L_a+zL_\beta+L_d \right),\\
L_e={x}^{2} \left( {z}^{2}+zL_\gamma+L_e   \right) ,\\
L_f={x}^{2} \left( zL_a+L_f \right) , 
\end{array}$$
and the \gf\ of completed polygons is
$$
S_3(x,y)={z}^{3}L_a+{z}^{2}L_\beta+zL_\gamma+z  L_d  +{z}^{2}L_f 
=\frac{xyN(x,y)}{D(x,y)}
$$
where
\begin{multline*}
N(x,y)=
3 \left( x+1 \right) ^{2} \left( 1-x \right) ^{5}
+ \left( 5\,x+2 \right)  \left( 2\,x-1 \right)  \left( x+1 \right)
^{2} \left(  x-1 \right) ^{3}{y}\\
-\left( x-1 \right)  \left( 6\,{x}^{6}+4\,{x}^{5}
-18\,{x}^{4}-6\,{x}^{3}+11\,{x}^{2}+8\,x+1 \right) {y}^{2}\\
-{x} \left( x+1 \right)  \left( 2\,{x}^{4}+6\,{x}^{3}-8\,{x}^{2}+4
\,x+1 \right) {y}^{3}
\end{multline*}
and
\begin{multline*}
D(x,y)=
\left( x+1 \right) ^{2}  \left( x-1 \right) ^{6}
-x \left( 1+4\,x \right)  \left( x+1 \right) ^{2} \left( x-1 \right) ^{4}y\\
+{x}^{2} \left( 3\,{x}^{4}+4\,{x}^{3}-6\,{x}^{2}-8\,x-3 \right) 
 \left( x-1 \right) ^{2}{y}^{2}\\
+{x}^{3} \left( x+1 \right)  \left( {x}^{3}+3\,{x}^{2}-5\,x+3 \right) {y}^{3}
.
\end{multline*}
By setting $x=y=t$, we obtain the half-perimeter \gf\ of SAP in a
3-strip,
$$
S_3(t)={\frac 
{{t}^{2} \left( -8\,{t}^{9}+4\,{t}^{8}+10\,{t}^{7}-20\,{t}^{6}-
{t}^{5}-{t}^{4}+7\,{t}^{3}+3\,{t}^{2}-7\,t+3 \right) }{4\,{t}^{10}-2\,
{t}^{9}-5\,{t}^{8}+8\,{t}^{7}-{t}^{6}+2\,{t}^{5}-4\,{t}^{4}+2\,{t}^{3}
+3\,{t}^{2}-4\,t+1}}
$$
 and, by looking at the smallest pole of this series, we also obtain
 the (very weak) 
 lower bound   $1.68...$ on the growth constant of square lattice
 self-avoiding polygons. 

The above method has been automated by Zeilberger~\cite{zeilberger-skinny}.
It is not hard to see that the number of states required to count
polygons in a $k$-strip grows like $3^k$, up to a power of $k$. This
prevents one from applying this method for large values of $k$. Better
bounds for growth constants may be obtained via the \emm finite
lattice method, 
described in Chapter~7. A further improvement is
obtained by looking at a cylinder rather than a strip~\cite{rote}.

%===================================================
\subsection{$q$-Analogues}
\label{sec:analogue-linear}
%===================================================
By looking at the height of the rightmost column of Ferrers diagrams,
we have described a 
linear construction of these polygons that proves the rationality of their perimeter
g.f. (Fig.~\ref{bb:fig:ferrers-perimeter}). 
 Let us  examine what happens when we try to
keep track of the area  in this construction.

They key point is that the area  \emm increases by the width of the \p, when
we duplicate the bottom row. (In contrast, the half-perimeter
simply increases by 1 during this operation.) This observation gives
the following functional equation for the complete g.f. of Ferrers diagrams:
$$
F(x,y,q)= xyq+ xqF(x,y,q)+  y F(xq,y,q).
$$
This is a $q$-\emm analogue, of the equation defining
$F(x,y,1)$, derived from~\eqref{ferrers-anisotropic}. This equation is
no longer linear, but it can be solved easily by iteration: 
\begin{eqnarray}
  F(x,y,q)&= &\frac{xyq}{1-xq} + \frac y {1-xq} F(xq,y,q)\nonumber\\
&= &\frac{xyq}{1-xq} + \frac y {1-xq}\frac{xyq^2}{1-xq^2}+ \frac y
  {1-xq}  \frac y {1-xq^2}F(xq^2,y,q)\label{bb:ferrers-complete} \\
&=& \sum_{n\ge 1} \frac{xy^n q^n}{(xq)_n}\nonumber
\end{eqnarray}
with 
$$
(xq)_0=1 \quad \hbox{ and }  \quad (xq)_n=(1-xq) (1-xq^2) \cdots (1-xq^n).
$$
Similarly, for the stack polygons of
Fig.~\ref{bb:fig:stacks-perimeter}, one obtains:
$$
\begin{array}{lcl}
S(x,y,q)&=&xyq+xqS(x,y,q)+S_+(x,y,q),\\
 S_+(x,y,q)&=& y S(xq,y,q)+xqS_+(x,y,q).
\end{array}
$$
Eliminating the series $S_+$ gives
\begin{eqnarray*}
S(x,y,q)&=& \frac{xyq}{1-xq} + \frac y {(1-xq)^2} S(xq,y,q)\\
&=& \sum_{n\ge 1} \frac{xy^n q^n}{(xq)_{n-1}(xq)_n}.
\end{eqnarray*}

In Section~\ref{bb:sec:layer} we present a  systematic approach for
counting classes of \cc\ \ps\ by perimeter and area.

%
%%%%%%%%%%%%%%%%%%%%%%%%%%%%%%%%%%%%%%%%%%%%%%%%%%%%%%%%%%
\section{Algebraic models and algebraic series}
\label{sec:alg}
%%%%%%%%%%%%%%%%%%%%%%%%%%%%%%%%%%%%%%%%%%%%%%%%%%%%%%%%%%

\subsection{A basic example: bargraphs counted by perimeter}
\label{sec:bargraphs-perimeter}
\begin{figure}[htb]
\begin{center}
 \scalebox{0.8}{\input{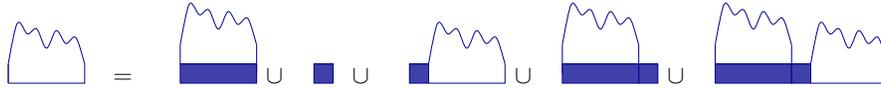}}
\end{center}
\caption{A second recursive construction of bargraphs.}
\label{bb:fig:bg-perimeter}
\end{figure}

Let us return to \bgs. The linear description used in
Section~\ref{sec:bg-area} to count them by  area cannot be directly
recycled to count them by perimeter: indeed, when we add a cell at
the top of the last column, how do we know if we increase the
perimeter, or not? Instead, we are going to scan the \p\ from left to
right, and \emm factor, it into two smaller bargraphs as soon as we meet a
column of height 1 (if any). If there is no such column, deleting the
bottom row of the polygon leaves another bargraph. This description is
schematized in Fig~\ref{bb:fig:bg-perimeter}.

Let $\cB$ be the set of words  over the alphabet $\{n,s,e\}$ that
naturally encode the top boundary of
bargraphs, from the SW to the SE corner. 
Fig.~\ref{bb:fig:bg-perimeter} translates into the
following recursive description, where the unions are disjoint: 
\begin{equation}\label{rec-bg-perimeter}
\cB= n \cL s \quad \hbox {with } \quad 
 \cL=  n \cL s \cup \{e \} \cup e \cL \cup n \cL s e \cup n \cL s e \cL.
\end{equation}
This implies that the anisotropic perimeter
 g.f. of \bgs\ satisfies
$$
\left\{
\begin{array}{lll}
 B(x,y)&=& y L(x,y),\\
L(x,y)&=& y L(x,y) + x+ x L(x,y) + xy L(x,y)  + xy L(x,y)^2.
\end{array}\right.
$$
These equations are readily solved and yield:
\begin{equation}\label{bb:bg-perimeter-sol}
B(x,y)= \frac{1-x-y-xy-\sqrt{(1-y)((1-x)^2-y(1+x)^2)}}{2x}.
\end{equation}
Thus the perimeter g.f. of \bgs \ is algebraic, and its algebraicity
is explained combinatorially by the recursive description of
Fig.~\ref{bb:fig:bg-perimeter}. 

%====================================================
Note that one can  directly translate this description into an
algebraic equation satisfied by $B(x,y)$, without using the language
$\cB$. This language is largely a convenient tool to highlight the
\emm algebraic structure, of bargraphs. The translation of
Fig.~\ref{bb:fig:bg-perimeter} into an equation proceeds as follows: there
are two types of bargraphs, those that have at least one column of
height 1, and the others, which we call \emm thick, bargraphs. Thick
bargraphs are obtained by duplicating the bottom row of a
general bargraph, and are thus counted by $yB(x,y)$. Among bargraphs
having a column of height 1, we find the single cell bargraph
(g.f. $xy$), and then those of width at least 2. The latter class can
be split into 3 disjoint classes:
\begin{itemize}
\item [--] the first column has height 1: these bargraphs are
  obtained by adding a cell to the left of any general bargraph, and
  are thus counted by $xB(x,y)$,
\item[--] the last column is the only column of height 1; these
  bargraphs are obtained  by adding a cell to the right of a thick
  bargraph, and   are thus counted by $xyB(x,y)$,
\item[--] the first column of height 1 is neither the first column,
  nor the last column. Such bargraphs are obtained by concatenating a
  thick bargraph, a cell, and a general bargraph; they are counted by
  $xB(x,y)^2$. 
\end{itemize}
This discussion directly results  in the equation
\begin{equation}
\label{bargraph-eq}
B(x,y)=yB(x,y)+xy+xB(x,y)+xy B(x,y)+ xB(x,y)^2.
\end{equation}
%====================================================

%================================================
\subsection{Algebraic  objects}
%================================================

The above description  of bargraphs  involved two constructions: 
\begin{enumerate}
  \item taking disjoint unions of sets,
\item taking cartesian products of sets.
\end{enumerate}
For two classes $\cA_1$ and $\cA_2$, the element $(a_1,a_2)$ of the
product  $\cA_1\times\cA_2$ is seen as the  concatenation of the
objects $a_1$ and $a_2$. 
We will say that a class of objects is \emm algebraic, if it admits a
non-ambiguous recursive description based on disjoint unions and
cartesian products. It is  assumed that the size of the objects is
\emm additive, for the concatenation. 
For instance,~\eqref{rec-bg-perimeter} gives an
algebraic description of the words of $\cL$ and $\cB$.

 In the case of linear
constructions, the only concatenations that were allowed were between
one object and a single atom. As we can now concatenate two objects,
algebraic constructions generalize linear
constructions.  In terms of g.f.s, concatenating objects of two
classes means taking the product of the corresponding g.f.s. Hence
the g.f. of an algebraic class will always be the first
component of  the
solution of a polynomial system of the form:
$$
A_i=P_i(t,A_1, \ldots, A_k) \quad \hbox{ for } \quad 1\le i \le k,
$$
where $P_i$ is a polynomial with coefficients in $\ns$. 
Such series are called $\ns$-\emm algebraic,, and are closely related to the theory
of \emm context-free languages,. We refer to~\cite{salomaa} for details on
these languages, and to~\cite{mbm-icm} for a discussion of
$\ns$-algebraic series in enumeration.

%================================================
\subsection{More algebraic models}
%================================================
\label{sec:algebraic-more}

In this section we present three problems that can be solved via an
algebraic decomposition: staircase polygons, then column-convex
polygons counted by perimeter (and the subclass of directed
column-convex polygons), and finally directed polyominoes counted
by area. 

%
%================================================
\subsubsection{Staircase polygons by perimeter}
\begin{figure}[htb]
\begin{center}
 \scalebox{0.8}{\input{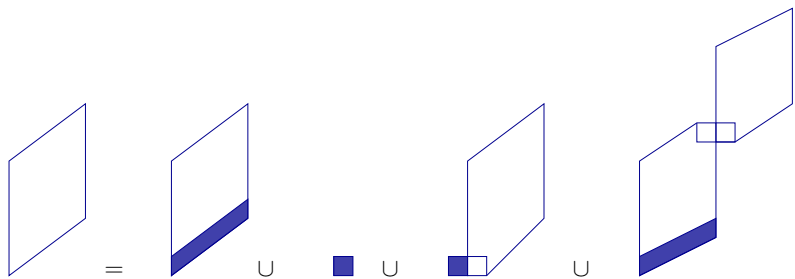}}
\end{center}
\caption{A recursive construction of staircase polygons.}
\label{bb:fig:staircase-perimeter}
\end{figure}

%================================================
In  \secref{sec:classes} we defined 
staircase  polygons  through their directed and convexity
properties.
See Fig.~\ref{bb:fig:zoo}(b) for an example. We describe here a recursive
construction of these polygons, illustrated in
Fig.~\ref{bb:fig:staircase-perimeter}. It is 
analogous to the construction of bargraphs described at the end of
Section~\ref{sec:bargraphs-perimeter} and illustrated in
Fig.~\ref{bb:fig:bg-perimeter}. Denote by $S(x,y)$ the anisotropic
perimeter \gf\ of staircase polygons.

We say that a staircase polygon is \emm thick, if deleting the bottom
cell of each column gives a staircase polygon of the same width. These
thick polygons are obtained by duplicating the bottom
cell in each column of a staircase polygon, so that their \gf\ is
$yS(x,y)$.

Among  non-thick staircase polygons, we find the single cell polygon
(g.f. $xy$),
and then those of width at least 2. Let $P$ be in the latter class,
and denote its columns $C_1, \ldots , C_k$, from left to right. The
fact that $P$ is not thick means that there exist two consecutive
columns, $C_i$ and $C_{i+1}$, that overlap by one edge only. Let
$i$ be minimal for this property. Two cases occur:
\begin{enumerate}
\item [--] the first column has height 1. In particular,  $i=1$.
  These polygons are obtained 
  by adding a cell to the bottom left of any general staircase
  polygon, and are thus counted by $xS(x,y)$,
\item [--] otherwise, the columns $C_1, \ldots, C_i$ form a thick staircase
  polygon, and $C_{i+1}, \ldots, C_k$ form a general staircase
  polygon. Concatenating these two polygons in such a way that they
  share only one edge gives the original polygon $P$. Hence the
  g.f. for this case is $S(x,y)^2$.
\end{enumerate}
This discussion gives the equation
$$
	S(x,y)=yS(x,y)+ xy+ x S(x,y) + S(x,y)^2\\
$$
so that
\begin{eqnarray*}
S(x,y)&=&\frac{1}{2} \left(1-x-y-\sqrt{1-2 x-2 y-2  xy+x^2+y^2 }
\right)  \\
&=&\sum_{p,q\ge 1} \frac 1 {p+q-1} {{p+q-1} \choose p} {{p+q-1}
  \choose q} x^p y^q.
\end{eqnarray*}
This expansion can be obtained using the Lagrange inversion formula~\cite{mbm:convex-languages}.
The isotropic semi-perimeter g.f. is obtained by setting $t=x=y$:
$$
S(t,t)=\frac{1}{2} \left(1-2 t-\sqrt{1-4 t} \right)
=\sum_{n\ge 1} C_n t^{n+1}
$$
where $C_n={{2n}\choose n}/(n+1)$ is the $n^{\hbox{th}}$ Catalan number.
The same approach can be applied to more general classes of convex
polygons, like directed-convex polygons and general convex
polygons. See for instance~\cite{mbm:convex-languages,delest-viennot}.

%
%================================================
\subsubsection{Column-convex polygons by perimeter}
\begin{figure}[htb]
\begin{center}
 \scalebox{0.8}{\input{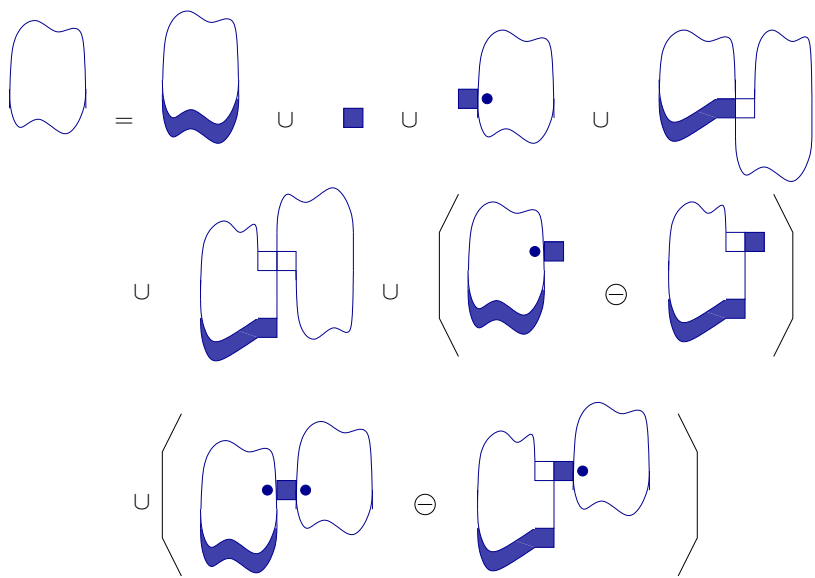}}
\end{center}
\caption{A recursive construction of column-convex polygons.}
\label{bb:fig:cc-perimeter}
\end{figure}

%================================================
We now apply a similar treatment to the perimeter enumeration of 
column-convex polygons (cc-polygons for short).
%See  Fig.~\ref{bb:fig:poly-ani}(c) for an example. 
Their area g.f. was found in Section~\ref{sec:linear-more}. 
Let $\cC$ denote the set of these
polygons, and $C(x,y)$ their anisotropic perimeter \gf. Our recursive
construction  requires us  to introduce two
additional classes of polygons. The first one, $\cC_1$, is the set of
cc-polygons  in which one cell of the last column is marked.  The
corresponding g.f. is denoted $C_1(x,y)$. Note that, by symmetry, this
series also counts cc-polygons where one cell of the \emm first, column
is marked. Then, $\cC_2$ denotes the set of
cc-polygons  in which one cell of the first column is marked (say, with
a dot),  and one cell of the last column is marked as well (say, with a
cross).
The corresponding g.f. is denoted $C_2(x,y)$.  Our recursive
construction of the polygons of $\cC$ is 
illustrated in
Fig.~\ref{bb:fig:cc-perimeter}. 

We say that a cc-polygon is \emm thick, if deleting the bottom
cell of each column gives a cc-polygon of the same width. These
thick polygons are obtained by duplicating the bottom
cell in each column of a cc-polygon, so that their \gf\ is
$yC(x,y)$.

Among  non-thick cc-polygons, we find the single cell polygon
(g.f. $xy$),
and then those of width at least 2. Let $P$ be in the latter class,
and denote its columns $C_1, \ldots , C_k$, from left to right. The
fact that $P$ is not thick means that there exist two consecutive
columns $C_i$ and $C_{i+1}$ that overlap by one edge only. Let
$i$ be minimal for this property. Two cases occur:
\begin{enumerate}
\item [--] the first column has height 1. In particular,  $i=1$.
  These polygons are obtained 
  by adding a cell to the  left of any cc-polygon having a
  marked cell in its first column, next to the marked cell. They are
  thus counted by $xC_1(x,y)$, 
\item [--] otherwise, the columns $C_1, \ldots, C_i$ form a thick 
  cc-polygon $P_1$, and the columns $C_{i+1}, \ldots, C_k$ form a general
  cc-polygon $P_2$.
  There are several ways of concatenating these two polygons
 in such a way they   share only one edge:
 \begin{enumerate}
 \item  [--] either the shared edge is at the bottom of $C_i$ and at the
   top of $C_{i+1}$: such polygons are counted by $C(x,y)^2$,
 \item  [--] or the shared edge is at the top of $C_i$ and at the
   bottom of $C_{i+1}$: such polygons are also counted by $C(x,y)^2$,
\item  [--] if $C_{i+1}$ has height at least 2, there are no other
  possibilities. However, if $C_{i+1}$ consists of one cell only, this
  cell may be adjacent to \emm any, cell of $C_i$, not only to the bottom
  or top ones. The case where $C_{i+1}$ is the last column of $P$ is
  counted by $xy(C_1(x,y)-C(x,y))$. The case where $i+1<k$ is counted
  by $x (C_1(x,y)^2-C(x,y)C_1(x,y))$.
 \end{enumerate}
\end{enumerate}
Let us drop the variables $x$ and $y$ in the series $C$, $C_1$ and $C_2$. The above
discussion gives the equation:
$$
C= yC+xy+xC_1+2C^2+ xy(C_1-C)+x (C_1^2-CC_1).
$$
The construction of Fig.~\ref{bb:fig:cc-perimeter} can now be recycled to obtain an
equation for the series $C_1$, counting cc-polygons with a marked cell
in the last column. Note that the first case of the figure (thick
polygons) gives 
rise to two terms, depending on whether the marked cell is one of the
duplicated cells, or not:
$$
C_1= y(C+C_1)+xy+xC_2+2CC_1+ xy(C_1-C)+x (C_1C_2-CC_2).
$$
We need a third equation, as three series (namely $C$, $C_1$ and $C_2$) are now
involved. There are two ways to obtain a third equation:
\begin{enumerate}
\item [--] either we interpret $C_1$ as the g.f. of cc-polygons where
  one cell is marked in the \emm first, column. The construction of
  Fig.~\ref{bb:fig:cc-perimeter} gives:
  \begin{multline*}
    C_1= y(C+C_1)+xy+xC_1+2(C+C_1)C+ xy\big((C_1-C)+(C_2-C_1)\big)\\
+x \big( (C_1^2-CC_1)+ (C_2C_1-C_1^2)\big).
  \end{multline*}
Note that now many cases give rise to two terms in the equation.
\item [--] or we work out an equation for $C_2$ using the
  decomposition of Fig.~\ref{bb:fig:cc-perimeter}. Again, many of the
  cases schematized in this figure give rise to several terms. In
  particular, the first case (thick polygons) gives rise to 4 terms:
  \begin{multline*}
    C_2= y(C+2C_1+C_2)+xy+xC_2+2(C+C_1)C_1+ xy\big((C_1-C)+(C_2-C_1)\big)\\
+x \big( (C_1C_2-CC_2)+ (C_2^2-C_1C_2)\big).
  \end{multline*}
\end{enumerate}
Both strategies  of course  give the same equation for $C\equiv C(x,y)$,
after the elimination of $C_1$ and $C_2$:
\begin{multline*}
\left( -5\,xy-18+2\,x{y}^{2}-18\,{y}^{2}+36\,y+2\,x \right) {C}^{4}\\
+ \left( y-1 \right)  \left( 5\,x{y}^{2}-21\,{y}^{2}+42\,y-14\,xy+5\,x-
21 \right) {C}^{3}\\
+2\, \left( y-1 \right) ^{2} \left( -4\,{y}^{2}+2\,x
{y}^{2}+8\,y-7\,xy-4+2\,x \right) {C}^{2}\\
+ \left( y-1 \right) ^{3}
 \left( x{y}^{2}-{y}^{2}+2\,y-6\,xy+x-1 \right) C-xy \left( y-1
 \right) ^{4}=0.
\end{multline*}
This quartic has 4 roots, among which the g.f. of cc-polygons can be
identified by checking the first few coefficients. This series turns out
to be unexpectedly simple: 
$$
C(x,y)=\left( 1-y \right)  \left( 1-
\frac{2\,\sqrt {2}}{ 3\,\sqrt {2}-\sqrt {1+x+\sqrt { \left( 1-x \right) ^{2}-16\,{\frac {xy}{ \left( 1-y
 \right) ^{2}}}}} }\right).
$$
Fereti\'c has provided direct combinatorial explanations
for this formula~\cite{feretic-2,feretic}. The algebraic equation satisfied by $C(t,t)$ was first\footnote{Eq.~(32) in~\cite{delest88} has an
error: the coefficient of $t^5c^3$ in $p_2$
should be $-40$ instead of  $+40$.}  obtained (via
a context-free language)
in~\cite{delest88}.  The method we have used is detailed in~\cite{duchi-rinaldi}.
%  The first few terms of the perimeter \gf\ are
%
%$$	  C(t,t)=t^2+2t^3+7t^4+28t^5+122t^6+558t^7+2641t^8+12822t^9+ O(t^{10}).$$

% 

%================================================
\subsubsection{Directed column-convex polygons by perimeter}
%================================================

\begin{figure}[htb]
\begin{center}
 \scalebox{0.8}{\input{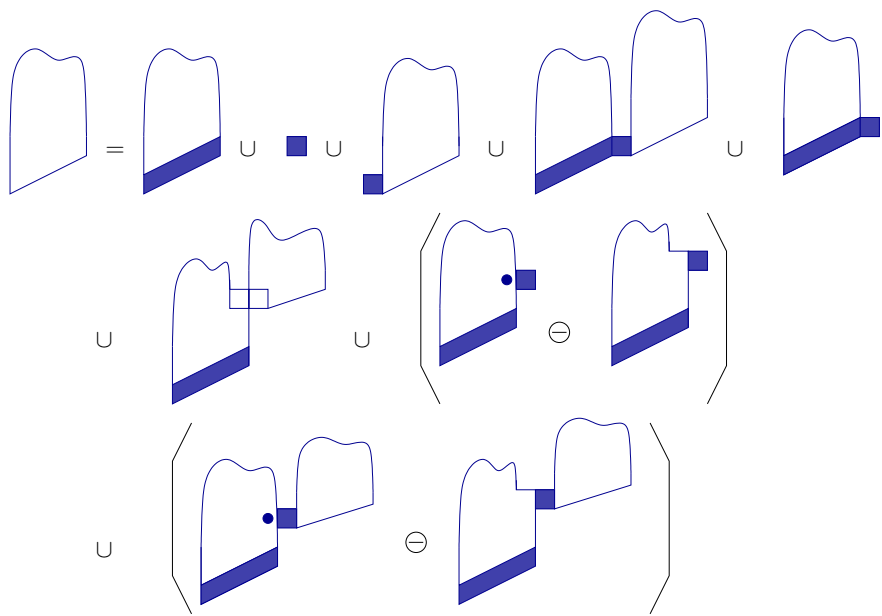}}
\end{center}
\caption{A recursive construction of directed column-convex polygons.}
\label{bb:fig:dc-perimeter}
\end{figure}

It is not hard to restrict the construction of Fig.~\ref{bb:fig:cc-perimeter} to \emm
directed, cc-polygons.  This is illustrated in
Fig.~\ref{bb:fig:dc-perimeter}. Note that the case where the
columns $C_i$ and $C_{i+1}$ share the bottom edge of $C_i$ (the fourth
case in Fig.~\ref{bb:fig:cc-perimeter}) is only possible if $C_{i+1}$ has height 1. 
Moreover,  only one additional series is needed,
namely that of directed cc-polygons marked in the last column
($D_1$). 

 One obtains the following equations:
 \begin{eqnarray*}
D&=&
yD+xy+xD+xD^2+xyD+D^2+xy(D_1-D)+x(D_1-D)D,
\\
D_1&=&
y(D+D_1)+xy+xD_1+xDD_1+xyD+DD_1+xy(D_1-D)\\
&&\hskip 80mm +x(D_1-D)D_1.
\end{eqnarray*}
Eliminating $D_1$ gives a cubic equation for the series $D\equiv D(x,y)$:
$$
{D}^{3}+ 2\left( y-1 \right) {D}^{2}+ \left( y-1 \right)  \left( x+y
-1 \right)  D  +xy \left( y-1 \right) =0.
$$
%\begin{align*}\label{eq:dccpgf}	
%	&P_3\, \dccgf^3+P_2\, \dccgf^2+P_1\, \dccgf+P_0\, =\, 0\\
%	\intertext{where}
%	P_0&= x y^3+x^3 y-x y\\
%	P_1&=y^2 x^4+x^4+2 y^2 x^2-2 x^2+y^4-2 y^2+1\\
%	P_2&=y^2 x^4+y^3 x^3+2 y x^3-y^2 x^2+2 y^3 x-2 y x\\
%	P_3&=x^3 y^3+x^2 y^2
%\end{align*}
This equation was first obtained in~\cite{delest:1993bh}.
The first few terms of the semi-perimeter generating function are
\begin{equation*} 
	D(t,t)=t^2+2 t^3+6 t^4+20t^5+71 t^{6}+263 t^{7}+1005 t^{8}+3933 t^{9}  +\cdots
\end{equation*}
%1, 2, 6, 20, 71, 263, 1005, 3933, 15684
% Somewhat surprisingly   the coupled equations \eqref{eq:dccpgf} can be solved using a multivariate  form of Lagrange inversion due to Good \cite{good:1960lr} to give the number of directed column-convex polygons with $k$ columns and perimeter $n$, $p_{n,k}$ as
%
%\begin{equation*}%\label{eq:dcccoeef}
%	p_{n,k}=\frac{1}{k}\sum_{r\ge0} \binom{2k}{k-r-1}\binom{2k+r-1}{r}\binom{n-2}{k+r-1}
%\end{equation*}
%

%================================================
\subsubsection{Directed polyominoes by area}
%================================================

Let us move to a
 class that admits a  neat, but non-obvious,
algebraic structure:  directed \pols\ counted by area. 
This structure was discovered when Viennot developed  the theory of
{\em heaps}~\cite{viennot-heaps}. 
 Intuitively, a heap is obtained by dropping vertically some
solid pieces,  one after the other. Thus, a piece lies either on
the ``floor'' (when it is said to be {\em minimal}), or at
least partially covers another piece.

\begin{figure}[hbt] \begin{center} 
 \scalebox{0.8}{\input{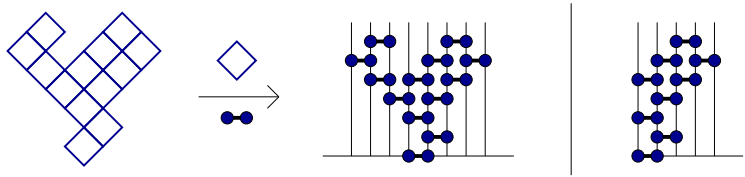}}
\end{center} 
\caption{\emm Left,: A directed \pol\ and the associated pyramid. 
\emm Right,: a half-pyramid.}
\label{bb:fig:dimers} 
\end{figure}

Directed \pols\ {\em are}, in essence, heaps. To see this,
replace every cell of the \pol\
by a  {\em dimer}, after a 45 degree rotation (Fig.~\ref{bb:fig:dimers}). 
This gives a  heap with  a unique minimal piece. 
Such heaps are called {\em pyramids}. 
If the columns to the left of the
minimal piece contain no dimer, we say we have a \emm half-pyramid, 
(Fig.~\ref{bb:fig:dimers}, right).

\begin{figure}[hbt] \begin{center} 
 \scalebox{0.8}{\input{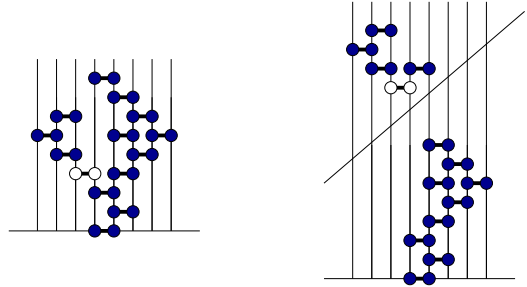}}
\end{center} 
\caption{A factorization of a pyramid into a pyramid and a
  half-pyramid.
Observe that the highest dimer of the pyramid moves up as we lift the
white dimer.}
\label{bb:fig-factorisation-heaps}
\end{figure}

The interest in  heaps lies in the existence of  a \emm product, of heaps: 
The product of two heaps is obtained by putting one heap above the
other and dropping its pieces.  Conversely, one can factor a heap by
pushing upwards one or several pieces. See an example in Fig.~\ref{bb:fig-factorisation-heaps}.
This product is the key in our  algebraic
description of directed \pols, or, equivalently, of pyramids of
dimers, as we now explain.

\begin{figure}[hbt] \begin{center} 
\scalebox{0.8}{\input{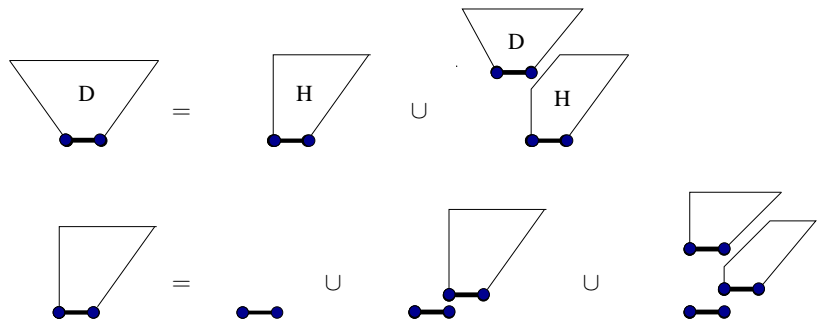}}
\caption{Decomposition of pyramids (top) and half-pyramids (bottom).}
\label{bb:fig:pyramide} 
\end{center} \end{figure}

A pyramid is either a   half-pyramid, 
or the product of a half-pyramid and  a pyramid
(Fig.~\ref{bb:fig:pyramide}, top). Let $D(q)$ denote the g.f.  of pyramids counted by the number of dimers, and
$H(q)$ denote the g.f. of half-pyramids.  Then $D(q)=H(q)(1+D(q))$. 

Now, a half-pyramid can be a single 
dimer. If it has several dimers, it is the product of a single dimer
and of one or two 
half-pyramids (Fig. \ref{bb:fig:pyramide}, bottom), which implies
$H(q)=q+qH(q)+qH^2(q)$.  

\noindent Note that $D(q)$ is also the area g.f. of directed \pols. A
straightforward computation  gives:
\begin{equation}\label{directed-area-sol}
D(q)=\frac{1}{2}\left( \sqrt{\frac{1+q}{1-3q}} -1 \right).
\end{equation}
This was first proved by Dhar~\cite{dhar-premiere}. The above proof is
adapted from~\cite{betrema-penaud}.
%================================================
\subsection{$q$-Analogues}
%================================================

By looking for the first column of height 1 in a bargraph, we have described an
algebraic construction of these \ps\ (Fig.~\ref{bb:fig:bg-perimeter})
that proves that their perimeter 
g.f. is algebraic (Section~\ref{sec:bargraphs-perimeter}). 
 Let us now examine what happens when we try to
keep track of the area of these \ps.

As in Section~\ref{sec:analogue-linear},
 the key observation is that the area behaves additively when one
concatenates two bargraphs, but \emm increases by the width of the \p, when
we duplicate the bottom row. (In contrast, the half-perimeter
simply increases by 1 during this operation.) This observation gives
rise to
the following functional equation for the complete g.f. of bargraphs:
\begin{multline}
B(x,y,q)= y B(xq,y,q) + xyq+ xqB(x,y,q)+ xyq B(xq,y,q) \label{q-alg-bg}\\
+ xq B(xq,y,q)B(x,y,q).
\end{multline}
This is a $q$-analogue of  Equation~\eqref{bargraph-eq} defining
$B(x,y,1)$. This equation is no longer algebraic, and it is not clear
how to solve it. It has been shown in~\cite{prellberg-brak} that it can
be linearized and solved using a certain Ansatz. We will show in
Section~\ref{bb:sec:bg-area-perimeter} a more systematic way to obtain
$B(x,y,q)$, which does not require any Ansatz.

%%%%%%%%%%%%%%%%%%%%%%%%%%%%%%%%%%%%%%%%%%%%%%%%%%%%%
\section{Adding a new layer: a versatile approach}
\label{bb:sec:layer}
%%%%%%%%%%%%%%%%%%%%%%%%%%%%%%%%%%%%%%%%%%%%%%%%%%%%%

In this section we describe a systematic construction that can be
used  to find the complete g.f. of 
many classes of \ps\ having
a convexity property~\cite{mbm-temperley}. The cost of this higher
generalization is twofold:
\begin{itemize}
  \item it is not always clear how to solve the functional equations
  obtained in this way,
\item in contrast with the constructions developed in
  Sections~\ref{sec:rat} and~\ref{sec:alg}, this approach does not
  provide combinatorial explanations for 
  the rationality/algebraicity of the corresponding g.f.s.
\end{itemize}
This type of construction is sometimes called \emm Temperley's approach,
since 
Temperley used it to write  functional equations for the generating
function of \cc \ \ps\ counted by
perimeter~\cite{temperley}. But it also occurs, in a more complicated form, in other ``old'' papers~\cite{bender,klarner-rivest}.  We would prefer to see a more precise
terminology, like \emm layered approach,.

%%%%%%%%%%%%%%%%%%%%%%%%%%%%%%%%%%%%%%%%%%%%%%%
\subsection{A basic example: bargraphs  by perimeter and area}
\label{bb:sec:bg-area-perimeter}
%%%%%%%%%%%%%%%%%%%%%%%%%%%%%%%%%%%%%%%%%%%%%%%

We return to our favourite example of bargraphs, and we now aim to 
find the complete g.f. $B(x,y,q)$ of this class of \ps. We have
just seen that the algebraic description of
Fig.~\ref{bb:fig:bg-perimeter} leads to the
$q$-algebraic equation~\eqref{q-alg-bg}, which is not obvious to
solve.
The linear description of Section~\ref{sec:bg-area} cannot be directly
exploited 
either: in order to decide whether the addition of a cell at the
top of the last column increases the perimeter or not, we need to
know which of the last two columns is higher.

We present here a variation of this linear construction that allows us
to count \bgs\ by area and perimeter, \emm provided we also take into
account the right height, by a new variable $s$. By \emm right height,, we mean the height of the rightmost column. 
The g.f. we are interested in is now
$$
B(x,y,q,s)=\sum_{h\ge 1} B_h(x,y,q) s^h,
$$
where $B_h(x,y,q)$ is the complete g.f. of \bgs\ of right height $h$. 

\begin{figure}[htb]
\begin{center}
 \scalebox{0.8}{\input{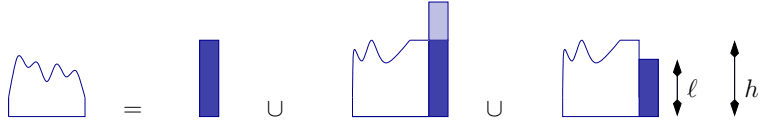}}
\end{center}
\caption{A third recursive construction of bargraphs.}
\label{bb:fig:bg-complete}
\end{figure}

Our new construction is illustrated in Fig.~\ref{bb:fig:bg-complete}. The
class $\cB$ of \bgs\ is split into three disjoint subsets:
\begin{enumerate}
  \item \bgs\ of width 1 (columns): the g.f. of this class is
  $xysq/(1-ysq)$,
\item \bgs\ in which the last column is at least as high as the
  next-to-last column. These \bgs\ are obtained by duplicating the
  last column of a \bg\ (which boils down to replacing $s$ by $sq$ in
  the series $B(x,y,q,s)$), and adding a (possibly empty) column at the
  top of the newly created column. The corresponding g.f. is thus
$$
\frac {x}{1-ysq} B(x,y,q,sq) .
$$
\item \bgs\ in which the last column is lower than the  next-to-last
  column. To obtain these, we start from a \bg, say of right height
  $h$, and add a new column of height $\ell<h$ to   the
  right. The g.f. of this third class is:
  \begin{eqnarray}
    x\sum_{h\ge 1} \left( B_h(x,y,q) \sum_{\ell=1}^{h-1}
    (sq)^{\ell}\right)\nonumber
&=&
x\sum_{h\ge 1}  \left( B_h(x,y,q) \ \frac{sq-(sq)^{h}}{1-sq}\right) \\
&=& x\ \frac{sq B(x,y,q,1) -  B(x,y,q,sq)}{1-sq}.  \label{geometric}
  \end{eqnarray}
\end{enumerate}
Writing $B(s)\equiv B(x,y,q,s)$, and putting together the three cases,
we obtain:
\begin{equation}\label{bb:eqfunc-bg}
B(s)= \frac {xysq}{1-ysq} +\frac {xsq}{1-sq}  B(1) + \frac
{xsq(y-1)}{(1-sq)(1-ysq)} B(sq). 
\end{equation}
This equation is solved in two steps: first, an iteration, similar to
what we did for Ferrers diagrams in~\eqref{bb:ferrers-complete}
(Section~\ref{sec:analogue-linear}), provides an expression for
$B(s)$ in terms of $B(1)$:
$$
B(s)= \sum_{n\ge 1} \frac{(xs(y-1))^{n-1} q^{n  \choose 2}}
{(sq)_{n-1}(ysq)_{n-1}} \left( \frac{xysq^n}{1-ysq^n}
+ \frac{xsq^n} {1-sq^n} B(1)\right).
$$
 Then, one sets $s=1$ to obtain the complete
 g.f. $B(1)\equiv B(x,y,q,1)$ of \bgs:
\begin{equation}\label{bb:bg-complete-sol}
B(x,y,q,1)= \frac {I_+}{1-I_-}
\end{equation}
with
$$
I_+= \sum_{n\ge 1} \frac{x^n(y-1)^{n-1} q^{{n+1}  \choose 2}}
{(q)_{n-1}(yq)_{n}} \quad\hbox{ and }  \quad
I_-=\sum_{n\ge 1} \frac{x^n(y-1)^{n-1} q^{{n+1}  \choose 2}}
{(q)_{n}(yq)_{n-1}} .
$$

%%%%%%%%%%%%%%%%%%%%%%%%%%%%%%%%%%%%%%%%%%%%%%%
\subsection{More examples}
\label{sec:newlayer-more}
%%%%%%%%%%%%%%%%%%%%%%%%%%%%%%%%%%%%%%%%%%%%%%%

In this section, we describe how to apply the layered
approach to two other classes of polygons: staircase and column-convex
polygons, counted by perimeters and area simultaneously. 
In passing we show how the  \emm difference, of g.f.s  resulting
from a geometric summation like~\eqref{geometric} can be explained
combinatorially by  an inclusion-exclusion argument. 

%================================================
\subsubsection{Staircase \ps}
%================================================

\begin{figure}[htbp] 
	\centering
		\includegraphics{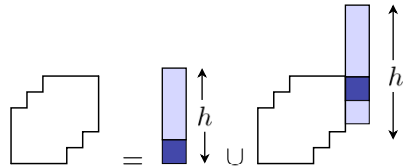}
	\caption{The two types of staircase polygons.}
	\label{fig_staircaseFTemperley}
\end{figure}

\begin{figure}[htbp] 
	\centering
		\includegraphics{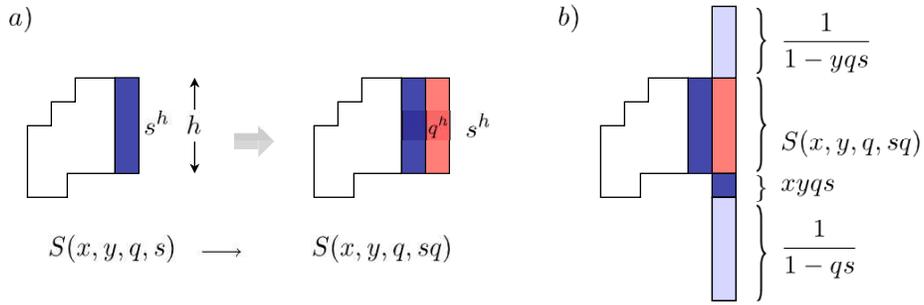}
	\caption{a) Replacing $s$ by $sq$ in $S(x,y,q,s)$ duplicates
          the last column of the polygon. b) Generating function 
of ``bad'' configurations.}
	\label{fig_staircaseFTemperley2}
\end{figure}

As with the bargraph example above, we define
an extended generating function which tracks the height of the
rightmost column of the staircase polygon, 
\begin{equation*}
	S(x,y,q,s)=\sum_{h\ge 1}S_h(x,y,q)s^h,
\end{equation*}
where $S_h(x,y,q)$ is the generating function of staircase polygons
with right height $h$. 
The set 
 of all staircase polygons can be partitioned into two
parts (Fig.~\ref{fig_staircaseFTemperley}):
\begin{enumerate}
	\item  those which have only one column: the g.f.\ for this
          class is $xyqs/(1-yqs)$,
	\item those which have more than one column: their g.f.\ 
 is obtained as the difference of the g.f.\ of two
          sets as follows. 
\end{enumerate}
Staircase polygons of width $\ell\ge 2$
can  be split into two objects: a staircase polygon formed of the
$\ell-1$ first columns, and the rightmost column.
 The left part has generating function $S(x,y,q,1)$ (ignoring the
 rightmost height), to which we then attach a column
 of cells. The attached column is
constrained in that it must not extend below the bottom of the
rightmost column of the left part. 
It is generated (see Fig.~\ref{fig_staircaseFTemperley}) by gluing a
descending column (with 
g.f.\ $1/(1-qs)$)  and an ascending column (with g.f.\ $1/(1-y qs)$)
to a single square (with g.f.\ $xyqs$). The single square is required
to ensure  that the column is not empty and is
 glued to the immediate right of the topmost square of the left part.
An important observation is that only the ascending column contributes
in the  increase of the vertical perimeter.
This gives the generating function 
\begin{equation*}
S(x,y,q,1)\cdot	xyqs\cdot \frac{1}{1-qs}\cdot \frac{1}{1-yqs}.
\end{equation*}
This construction 
however results in configurations which might have the rightmost
column extending below the rightmost column of the left part.
 We must thus subtract the contribution of these ``bad''
 configurations from the above g.f.
We claim that they  are generated by 
\[
S(x,y,q,sq)\cdot	xyqs\cdot \frac{1}{1-qs}\cdot \frac{1}{1-yqs}.
\]
The replacement of $s$ with $sq$ in $S(x,y,q,sq)$ 
is interpreted as adding a copy of the last
column of the left part,
as illustrated in Fig.~\ref{fig_staircaseFTemperley2}. 
The $xyqs$ factor is  interpreted as  attaching a new cell to the bottom of
the duplicated column (thus ensuring the rightmost column is strictly
below the rightmost column of the left part).
Finally, we add  a descending and an ascending column.  Again,
the height of the latter must not be taken into account in the
vertical perimeter.

Thus the final equation for the generating function is  
\begin{equation*}\label{eq_staircaseFunc}
	S(x,y,q,s)=\frac{xyqs}{1-yqs}+\bigl(S(x,y,q,1)-S(x,y,q,sq)\bigr) \frac{xyqs }{(1-qs)(1-yqs)}.
\end{equation*}
It can also be obtained via geometric sums, as was done for~\eqref{geometric}.
The equation is solved with the same two step process as for
bargraphs. First we iterate it 
to obtain $S(x,y,q,s)$ in terms of $S(x,y,q,1)$, 
and then we set $s$ to 1, obtaining
\begin{equation*}
	S(x,y,q,1)=y\frac{J_1}{J_0},
\end{equation*}
where $J_0$ and $J_1$ are  two $q-$Bessel functions \cite{gasper90}:
\begin{equation*}
	J_1(x,y,q)= \sum_{k\ge 1}(-1)^{k+1}\frac{x^k
          q^{\binom{k+1}{2}}}
{(q)_{k-1}(yq)_k}
\end{equation*}
and
\begin{equation*}
	J_0(x,y,q)= \sum_{k\ge 0}(-1)^{k}\frac{x^k q^{\binom{k+1}{2}}}
{(q)_{k}(yq)_k}.
\end{equation*}
Note, the appropriate limit as $q\to1$ leads to  standard Bessel
functions which  
are related to the generating function for
semi-continuous staircase polygons -- see~\cite{brak:1994qk} for details. 

%=================================================
\subsubsection{Column-convex \ps}
%================================================

The case of column-convex polygons is  more complex and
we will not give all the details but only discuss the primary
additional complication. We refer to~\cite{mbm-temperley} for a
complete solution. 
Like for staircase polygons,  a functional
equation for column-convex polygons can be obtained by considering
the rightmost (last) column.
 The position of the last column compared with the second last column  must be
carefully considered.  
Again there are several cases depending on whether  the top (resp.\
bottom) of the last column is strictly above, at the same level
 or below the
top (resp. bottom) of the second-last column. The case which leads to a type of
term that does not appear in the equation for staircase polygons is
the case where the 
top (resp. bottom) of the last column cannot be above (resp. below) the top
(resp. bottom) of  the second-last column. Thus we will only explain this %last
case which we will refer to as the \emm contained, case, as the last
column is somehow contained in the previous one.

If the generating function  for the column-convex polygons is
$C(s)=C(x,y,q,s)$ then we claim that 
the polygons falling into the contained case
are counted by the generating function 
\begin{equation}
	\frac{xsq}{1-sq}\frac{\partial C}{\partial
          s}(1)-\frac{xs^2q^2}{(1-sq)^2}\bigl(C(1)-C(sq)\bigr).
\label{eq-cc-temperley}
\end{equation}
Thus we see we now need a derivative
  of the generating function. 
As a polygon of right height $h$
  contributes $h$ times to the series ${\partial C}/{\partial s}(1)$,
   this series counts polygons with a marked cell in the
  rightmost column. 

\begin{figure}[htbp] 
	\centering
		\includegraphics{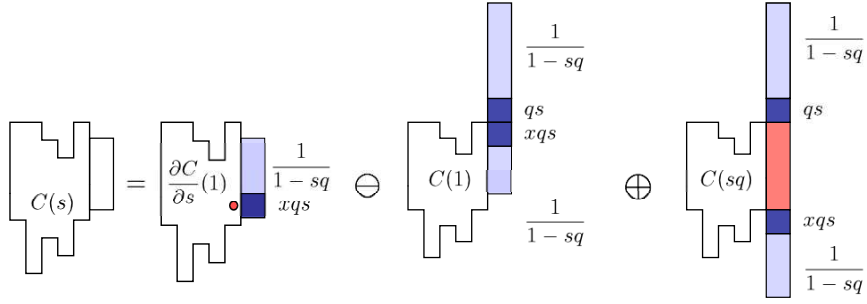}
	\caption{A schematic representation of the equation  for the case where the rightmost column does not extend above or below the second-last column.}
	\label{fig_ccFTemperleyEQ}
\end{figure}

Let us now explain this expression, which is illustrated
 in \figref{fig_ccFTemperleyEQ}. 
We consider a polygon as the concatenation of a left part
with a new (rightmost) column $C$. In the left part, we mark the cell of
the rightmost column that is at the same level as the bottom cell of
$C$.
So, starting from a marked polygon,
we first add a single square to the right of the marked cell
-- this gives a factor $xsq$. Above this square we then add an ascending 
column  which is generated by $1/(1-sq)$.
However, as with the staircase polygons, the resulting series
counts ``bad'' configurations, 
where the last column ends strictly higher than
the second last column.
 We subtract the contribution of these bad
configurations by generating them as shown on the second picture of
Fig.~\ref{fig_ccFTemperleyEQ}. This results in subtracting the term
$xq^2s^2C(1)/(1-sq)^2$. However, we have now subtracted too much! Indeed,
some configurations counted by the latter series have a 
rightmost column that ends below the second last column. We correct
this by adding the contribution of these configurations, which is
$xq^2s^2 C(sq)/(1-sq)^2$ (Fig.~\ref{fig_ccFTemperleyEQ}, right). This
establishes~\eqref{eq-cc-temperley} for the g.f. of the contained case.

The other cases are simpler, and in the same vein as what was needed for
staircase polygons. Considering all cases gives
\begin{multline}
	C(s)=\frac{xsyq}{1-syq}+
	\frac{xsq}{1-sq}\frac{\partial C}{\partial s}(1)
	+\frac{xs^2q^2(2y-syq-1)}{(1-sq)^2(1-syq)}C(1)\\
	+\frac{xs^2q^2(1-y)^2}{(1-sq)^2(1-syq)^2}C(sq).
\end{multline}
In order to solve this equation, we first iterate it to obtain $C(s)$
in terms of $C(1)$ and $C'(1) ={\partial C}/{\partial s}(1)$. 
Setting $s=1$ gives a  linear equation between $C(1)$ and
$C'(1)$. Setting $s=1$ \emm after having differentiated with respect
to $s$, gives a second linear equation between $C(1)$ and $C'(1)$. We
end up solving a linear system of size 2, and obtain $C(1)$ as a ratio
of two $2\times2$ determinants. The products of series that
appear in these determinants can be simplified, and the final expression reads
\begin{equation*}
C(x,y,q,1) =y\frac{  (1-y)X  }{  1+W+yX } 
\end{equation*}
where
\begin{align*}
 X&= \frac{xq}{(1-y)(1-yq)}
+\sum_{n\ge2}   \frac{    (-1)^{n+1}x^n(1-y)^{2n-4}q^{\binom{n+1}{2}} 
(y^2q)_{2n-2}    }
{ (q)_{n-1}\, (yq)_{n-2}\,   (yq)^2_{n-1}\,   (yq)_{n}\,   (y^2q)_{n-1}   }  
\intertext{and}
W&=\sum_{n\ge1}\frac{   (-1)^{n}
x^n(1-y)^{2n-3}q^{\binom{n+1}{2}} (y^2q)_{2n-1}   }
{   (q)_{n}\,  (yq)^3_{n-1}\,  (yq)_{n}\,  (y^2q)_{n-1}   }.
\end{align*}
The first solution, involving a more complicated expression, was given
in \cite{brak:1990vi}. The  one above
appears as Theorem 4.8 in~\cite{mbm-temperley}.

%================================================
\subsection{The kernel method}
%================================================

In Sections~\ref{sec:rat} and \ref{sec:alg}, we have  explained
combinatorially why the area g.f. of bargraphs, $B(1,1,q)$, and the
perimeter g.f. of bargraphs, $B(x,y,1)$, are respectively rational and
algebraic.  It is natural to examine whether  these properties can be
recovered from the construction of Fig.~\ref{bb:fig:bg-complete} and
the functional equation~\eqref{bb:eqfunc-bg}. 

As soon as we set $y=1$ in this equation, the main difficulty, that
is, the term $B(sq)$, disappears. We can then substitute $1$ for $s$ and
solve for $B(x,1,q,1)$, the width and area g.f. of bargraphs. This series is
found to be
$$
B(x,1,q,1)= \frac{xq}{1-q-xq}.
$$
From this, one also obtains a rational expression for the series
$B(x,1,q,s)$. The rationality of $B(x,1,q,1)$ also follows directly
from the expression~\eqref{bb:bg-complete-sol}: setting $y=1$ 
shrinks the series  $I_+$ and $I_-$ to simple rational functions.

How the \emm perimeter, g.f. of \bgs\ can be derived from the functional
equation~\eqref{bb:eqfunc-bg} is a more challenging question. Setting
$q=1$ gives
\begin{equation}\label{bb:eqfunc-bg-q=1}
B(s)= \frac {xys}{1-ys} +\frac {xs}{1-s}  B(1) + \frac
{xs(y-1)}{(1-s)(1-ys)} B(s). 
\end{equation}
This equation cannot be simply solved by setting $s=1$. Instead, the
solution uses the so-called \emm kernel method,, which has proved
useful in a rather large variety of enumerative problems in the past
10 years~\cite{hexacephale,banderier-flajolet,mbm-petkovsek,de-mier,feretic-svrtan,prodinger}. This method solves, in a systematic way, equations of
the form:
$$
K(s,x) A(s,x)= P(x,s,A_1(x), \ldots, A_k(x))
$$
where $K(s,x)$ is a polynomial in $s$ and the other indeterminates
$x=(x_1, \ldots, x_n)$, $P$ is a polynomial, $A(s,x)$ is an unknown
series in $s$ and the $x_i$'s, while the series $A_i(x)$ only depend on
the $x_i$'s. (It is assumed that the equation  uniquely defines all
these unknown series.) We refer to~\cite{mbm-petkovsek} for a general
presentation, and simply illustrate the method
on~\eqref{bb:eqfunc-bg-q=1}.
We group the terms involving $B(s)$, and multiply the equation by $(1-s)$
to obtain:
\begin{equation}\label{bb:eqfunc-bg-q=1-bis}
\left(1-s-\frac{xs(y-1)}{1-ys}\right) B(s)= \frac{xys(1-s)}{1-ys} + xsB(1).
\end{equation}
Let $S\equiv S(x,y)$ be the only \fps\ in $x$ and $y$ that satisfies
$$ S=1-\frac{xS(y-1)}{1-yS}.
$$
That is, 
$$
S=\frac{1-x+y+xy-\sqrt{(1-y)((1-x)^2-y(1+x)^2)}}{2y}.
$$
Replacing $s$ by $S$ in~\eqref{bb:eqfunc-bg-q=1-bis} gives an identity
between series in $x$ and $y$. By construction, the left-hand
side of this identity vanishes. This gives
\begin{multline*}
  B(1)\equiv B(x,y,1,1)= \frac{y(S-1)}{1-yS}\\
= \frac{1-x-y-xy-\sqrt{(1-y)((1-x)^2-y(1+x)^2)}}{2x},
\end{multline*}
and we have recovered the algebraic
expression~\eqref{bb:bg-perimeter-sol} of the perimeter g.f. of \bgs.

% 
%================================================
\section{Some open questions}
\label{sec:questions}
%================================================

We conclude this chapter with a list of open questions. As mentioned
in the introduction, the combination of convexity and direction
conditions 
gives rise to 35  classes of polyominoes, not all solved. But all
these classes are certainly not  equally interesting. The few problems
we present below have two important qualities: they do not seem
completely out of reach (we do not ask about the enumeration of all
polyominoes) and they have some special interest: 
they deal either with large classes of polyominoes, or with mysterious
classes (that have been solved in a non-combinatorial fashion), or
they seem to lie just at the border of what the available techniques
can achieve at the moment. 

%================================================
\subsection{The quasi-largest class of quasi-solved polyominoes}
%================================================

Let us recall that the growth constant of $n$-cell polyominoes is
conjectured to be a bit more than 4. More precisely, it is believed
that $p_n$, the number of such polyominoes, is equivalent to $\mu^n
n^{-1}$,
up to a multiplicative constant, with $\mu=4.06...$~\cite{jensen-guttmann-animals}.
  The techniques that provide lower bounds on $\mu$ involve looking
  at \emm bounded, polyominoes (for instance polyominoes lying in a strip of fixed
  height $k$) and a concatenation argument. See~\cite{rote} for a recent
  survey and the best published lower bound, $3.98...$.
 It is not hard to see that for $k$ fixed, 
these bounded polyominoes have a linear structure, and a rational \gf.  This
series is obtained either by adding recursively a whole ``layer'' to
the polyomino (as we did for self-avoiding polygons in
Section~\ref{sec:polygons-strip}), or by adding  
one cell at a time. The latter approach is 
usually more efficient (Chapter~6). 

What about  solved classes of polyominoes that do not depend on a
parameter $k$, and often have a more subtle structure? We have seen
in Section~\ref{sec:algebraic-more} that the g.f. of directed polyominoes is
algebraic, with growth constant  3. This is ``beaten'' by the growth
constant 3.20... derived from the rational g.f. of column-convex  polyominoes
(Section~\ref{sec:linear-more}). A 
generalization of directed polyominoes (called multi-directed
polyominoes) was introduced 
in~\cite{mbm-rechni-heaps} and proved to have a fairly complicated g.f.,
with growth constant about $3.58$. To our knowledge, this is the
largest growth constant reached from exact enumeration (again, apart
from the rational classes obtained by bounding column heights). 

However, in 1967, Klarner introduced a ``large'' class of polyominoes that
seems interesting and would  warrant a better understanding~\cite{klarner}. His 
definition is a bit unclear, and his solution is only partial, but the
estimate he obtains of the growth constant is definitely appealing:
about $3.72$. Let us mention that the triangular lattice version of
this mysterious class is solved in~\cite{mbm-rechni-heaps}. The growth
constant is found to be about $4.58$ (the growth constant of triangular lattice
animals is estimated to be about $5.18$, see~\cite{voege-guttmann}). 

%==============================================
\subsection{Partially directed polyominoes}
\begin{figure}[hbt]
\begin{center}
\includegraphics[width=12cm]{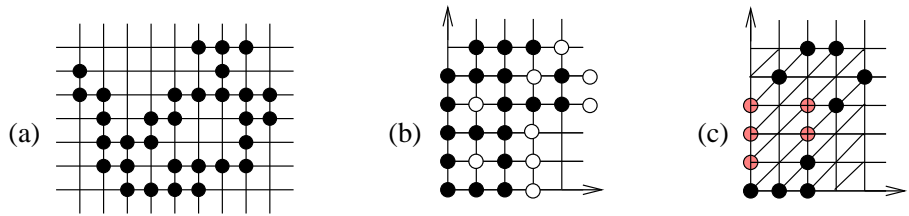}
\end{center}
\caption{(a) A partially directed animal. The source can be any point
  on the bottom row. (b) A directed animal on the
  square lattice, with the  right neighbours indicated in white. (c) A
  directed animal $A$ on the triangular lattice. The distinguished
  points are those having  (only) their South neighbour in $A$.}
\label{bb:fig:directed-variants}
\end{figure}

%==============================================
This is another generalization of directed polyominoes, with a very
natural definition:  the corresponding animal $A$ contains a source point
$v_0$ from which every other point can be reached by a path
formed of North, East and West steps, only visiting points of $A$
(Fig.~\ref{bb:fig:directed-variants}(a)). This model has a 
slight flavour of \emm heaps of pieces,, a notion that has already
proved useful in the solution of several polyomino
models (see Section~\ref{sec:algebraic-more}
and~\cite{betrema-penaud,mbm-viennot,mbm-rechni-heaps}). The growth
constant is estimated to be around 3.6, and, if proved, would thus improve that of
multi-directed animals~\cite{privman-barma}.

%================================================
\subsection{The right site-perimeter of directed animals}
%================================================

We wrote in the introduction that \emm almost all, solved classes of
polyominoes can be solved by one of the three main approaches we
present in this chapter. Here is one simple-looking result that we do
not know how to prove via these approaches (nor via any combinatorial
approach, to be honest). 

Take a directed animal $A$, and call a \emm neighbour, of $A$ any point
that does not lie in $A$, but could be added to $A$ to form a new
directed animal. The number of neighbours is the \emm site-perimeter,
of $A$. The \emm right site-perimeter, of $A$ is the number of
neighbours that lie one step to the right of a point of $A$. It was
proved in~\cite{mbm-dhar} that the g.f. of directed animals, counted by
 area and right site-perimeter, is a very simple extension
of~\eqref{directed-area-sol}:
$$
D(q,x)=\frac x 2 \left( \sqrt{\frac{(1+q)(1+q-qx)}{1-q(2+x)+q^2(1-x)}}-1\right).
$$
The proof is based on an equivalence with a one-dimensional gas model,
inspired by~\cite{dhar}.
It is easy to see that the right site-perimeter is also the number of
vertices $v$ of $A$ whose West neighbour is not in $A$. (By the West
neighbour, we mean the point at coordinates $(i-1,j)$   if
$v=(i,j)$). 

Described in these terms, this result has a remarkable counterpart for
 triangular lattice animals (Fig.~\ref{bb:fig:directed-variants}). Let us say
 that a point $(i,j)$ of the animal has a West (resp. South, South-West)
 neighbour in $A$ if the point $(i-1,j)$ (resp. $(i,j-1)$, $(i-1,j-1)$)
 is also in $A$.  Then the g.f. that counts these animals by the area
 and the number of points having a SW-neighbour (but no W- nor
 S-neighbour) is easy to obtain using  heaps of dimers and the ideas
 presented in Section~\ref{sec:algebraic-more}: 
$$
\tilde D(q,x)= \frac 1 2 \left( \sqrt{\frac{1+q-qx}{1-3q-qx}}\right).
$$
What is less easy, and is so far only proved via a correspondence
with a gas model, is that $\tilde D(q,x)$ also counts directed animals
(on the triangular lattice) by the area and the number of points
having a South neighbour (but no SW- nor W-neighbour). Any
combinatorial proof of this result would give a better understanding
of these objects. One possible starting point may be found in the
recent paper~\cite{borgne-marckert}, which sheds some combinatorial
light on the gas models involved in the proof of the above identities.

%==============================================
\subsection{Diagonally convex polyominoes} 
%==============================================
Let us conclude with a problem that seems to lie at the border of the
applicability of the third approach presented here (the layered
approach). 
In the enumeration of, say, column-convex polyominoes
(Section~\ref{sec:newlayer-more}), we have used the fact that deleting the last column
of such a polyomino gives another  \cc\ polyomino. This is no longer true
of a $d_-$-convex polyomino from which we would delete the last
diagonal (Fig.~\ref{bb:fig:poly-ani}(d)). Still, it seems that this class is
sufficiently well structured to be exactly enumerable. Note that this
difficulty vanishes when studying the restricted class of \emm
directed, diagonally convex 
polyominoes~\cite{mbm-habilitation,svrakic}, which behave
approximately like \cc \ polyominoes.

%%%%%%%%%%%%%%%%%%%%%%%%%%%%%%%%%%%%%%%%%%%%%%%%%%%%%%%%%%%%%%
\bibliographystyle{plain}

\bibliography{bibli-polygons}

\begin{thebibliography}{10}

\bibitem{alm-janson}
S.~E. Alm and S.~Janson.
\newblock Random self-avoiding walks on one-dimensional lattices.
\newblock {\em Comm. Statist. Stochastic Models}, 6(2):169--212, 1990.

\bibitem{andrews}
G.~E. Andrews.
\newblock {\em The theory of partitions}.
\newblock Addison-Wesley Publishing Co., Reading, Mass.-London-Amsterdam, 1976.
\newblock Encyclopedia of Mathematics and its Applications, Vol. 2.

\bibitem{hexacephale}
C.~Banderier, M.~Bousquet-M{\'e}lou, A.~Denise, P.~Flajolet, D.~Gardy, and
  D.~Gouyou-Beauchamps.
\newblock Generating functions for generating trees.
\newblock {\em Discrete Math.}, 246(1-3):29--55, 2002.

\bibitem{banderier-flajolet}
C.~Banderier and P.~Flajolet.
\newblock Basic analytic combinatorics of directed lattice paths.
\newblock {\em Theoret. Comput. Sci.}, 281(1-2):37--80, 2002.

\bibitem{rote}
G.~Barequet, M.~Moffie, A.~Rib{\'o}, and G.~Rote.
\newblock Counting polyominoes on twisted cylinders.
\newblock {\em Integers}, 6:A22, 37 pp. (electronic), 2006.

\bibitem{bender}
E.~A. Bender.
\newblock Convex {$n$}-ominoes.
\newblock {\em Discrete Math.}, 8:219--226, 1974.

\bibitem{betrema-penaud}
J.~B\'etr\'ema and J.-G. Penaud.
\newblock Mod\`eles avec particules dures, animaux dirig\'es et s\'eries en
  variables partiellement commutatives.
\newblock ArXiv:math.CO/0106210.

\bibitem{mbm-habilitation}
M.~Bousquet-M\'elou.
\newblock Rapport scientifique d'habilitation.
\newblock Report 1154-96, LaBRI, Universit\'e Bordeaux 1,
  {http://www.labri.fr/perso/lepine/Rapports\_internes}.

\bibitem{mbm:convex-languages}
M.~Bousquet-M{\'e}lou.
\newblock Codage des polyominos convexes et \'equations pour l'\'enum\'eration
  suivant l'aire.
\newblock {\em Discrete Appl. Math.}, 48(1):21--43, 1994.

\bibitem{mbm-temperley}
M.~Bousquet-M{\'e}lou.
\newblock A method for the enumeration of various classes of column-convex
  polygons.
\newblock {\em Discrete Math.}, 154(1-3):1--25, 1996.

\bibitem{mbm-dhar}
M.~Bousquet-M{\'e}lou.
\newblock New enumerative results on two-dimensional directed animals.
\newblock {\em Discrete Math.}, 180(1-3):73--106, 1998.

\bibitem{mbm-icm}
M.~Bousquet-M\'elou.
\newblock Rational and algebraic series in combinatorial enumeration.
\newblock In {\em Proceedings of the International Congress of Mathematicians},
  pages 789--826, Madrid, 2006. European Mathematical Society Publishing House.

\bibitem{mbm-guttmann}
M.~Bousquet-M{\'e}lou and A.~J. Guttmann.
\newblock Enumeration of three-dimensional convex polygons.
\newblock {\em Ann. Comb.}, 1(1):27--53, 1997.

\bibitem{mbm-petkovsek}
M.~Bousquet-M{\'e}lou and M.~Petkov{\v{s}}ek.
\newblock Linear recurrences with constant coefficients: the multivariate case.
\newblock {\em Discrete Math.}, 225(1-3):51--75, 2000.

\bibitem{mbm-rechni-heaps}
M.~Bousquet-M{\'e}lou and A.~Rechnitzer.
\newblock Lattice animals and heaps of dimers.
\newblock {\em Discrete Math.}, 258(1-3):235--274, 2002.

\bibitem{mbm-viennot}
M.~Bousquet-M{\'e}lou and X.~G. Viennot.
\newblock Empilements de segments et {$q$}-\'enum\'eration de polyominos
  convexes dirig\'es.
\newblock {\em J. Combin. Theory Ser. A}, 60(2):196--224, 1992.

\bibitem{brak:1990vi}
R.~Brak and A.~J. Guttmann.
\newblock Exact solution of the staircase and row-convex polygon perimeter and
  area generating function.
\newblock {\em J. Phys A: Math. Gen.}, 23(20):4581--4588, 1990.

\bibitem{brak:1994qk}
R.~Brak, A.~L. Owczarek, and T.~Prellberg.
\newblock Exact scaling behavior of partially convex vesicles.
\newblock {\em J. Stat. Phys.}, 76(5/6):1101--1128, 1994.

\bibitem{de-mier}
A.~de~Mier and M.~Noy.
\newblock A solution to the tennis ball problem.
\newblock {\em Theoret. Comput. Sci.}, 346(2-3):254--264, 2005.

\bibitem{italiens-conv-dir}
A.~Del~Lungo, M.~Mirolli, R.~Pinzani, and S.~Rinaldi.
\newblock A bijection for directed-convex polyominoes.
\newblock In {\em Discrete models: Combinatorics, Computation, and Geometry
  (Paris, 2001)}, {\em Discrete Math. Theor. Comput. Sci. Proc.}, pages
  133--144 (electronic). Maison Inform. Math. Discr., Paris, 2001.

\bibitem{delest:1993bh}
M.~Delest and S.~Dulucq.
\newblock Enumeration of directed column-convex animals with given perimeter
  and area.
\newblock {\em Croatica Chemica Acta}, 66(1):59--80, 1993.

\bibitem{delest88}
M.-P. Delest.
\newblock Generating functions for column-convex polyominoes.
\newblock {\em J. Combin. Theory Ser. A}, 48(1):12--31, 1988.

\bibitem{delest-viennot}
M.-P. Delest and G.~Viennot.
\newblock Algebraic languages and polyominoes enumeration.
\newblock {\em Theoret. Comput. Sci.}, 34(1-2):169--206, 1984.

\bibitem{dhar-premiere}
D.~Dhar.
\newblock Equivalence of the two-dimensional directed-site animal problem to
  {B}axter's hard square lattice gas model.
\newblock {\em Phys. Rev. Lett.}, 49:959--962, 1982.

\bibitem{dhar}
D.~Dhar.
\newblock Exact solution of a directed-site animals-enumeration problem in
  three dimensions.
\newblock {\em Phys. Rev. Lett.}, 51(10):853--856, 1983.

\bibitem{duchi-rinaldi}
E.~Duchi and S.~Rinaldi.
\newblock An object grammar for column-convex polyominoes.
\newblock {\em Ann. Comb.}, 8(1):27--36, 2004.

\bibitem{enting-guttmann-sap}
I.~G. Enting and A.~J. Guttmann.
\newblock On the area of square lattice polygons.
\newblock {\em J. Statist. Phys.}, 58(3-4):475--484, 1990.

\bibitem{feretic-2}
S.~Fereti{\'c}.
\newblock The column-convex polyominoes perimeter generating function for
  everybody.
\newblock {\em Croatica Chemica Acta}, 69(3):741--756, 1996.

\bibitem{feretic}
S.~Fereti{\'c}.
\newblock A new way of counting the column-convex polyominoes by perimeter.
\newblock {\em Discrete Math.}, 180(1-3):173--184, 1998.

\bibitem{feretic-festoon1}
S.~Fereti{\'c}.
\newblock An alternative method for {$q$}-counting directed column-convex
  polyominoes.
\newblock {\em Discrete Math.}, 210(1-3):55--70, 2000.

\bibitem{feretic-festoon2}
S.~Fereti{\'c}.
\newblock A {$q$}-enumeration of convex polyominoes by the festoon approach.
\newblock {\em Theoret. Comput. Sci.}, 319(1-3):333--356, 2004.

\bibitem{feretic-svrtan}
S.~Fereti{\'c} and D.~Svrtan.
\newblock On the number of column-convex polyominoes with given perimeter and
  number of columns.
\newblock In Barlotti, Delest, and Pinzani, editors, {\em Proceedings of the
  5th Conference on Formal Power Series and Algebraic Combinatorics (Florence,
  Italy)}, pages 201--214, 1993.

\bibitem{flaj-sedg2}
P.~Flajolet and R.~Sedgewick.
\newblock {\em Analytic Combinatorics}.
\newblock Preliminary version available at
  http://pauillac.inria.fr/algo/flajolet/Publications/books.html.

\bibitem{gasper90}
G.~Gasper and M.~Rahman.
\newblock {\em Basic hypergeometric series}, volume~35 of {\em Encyclopedia of
  Mathematics and its Applications}.
\newblock Cambridge University Press, Cambridge, 1990.

\bibitem{guttmann-prellberg}
A.~J. Guttmann and T.~Prellberg.
\newblock Staircase polygons, elliptic integrals, {H}eun functions and lattice
  {G}reen functions.
\newblock {\em Phys. Rev. E}, 47:R2233--R2236, 1993.

\bibitem{jensen-guttmann-sap}
I.~Jensen and A.~J. Guttmann.
\newblock Self-avoiding polygons on the square lattice.
\newblock {\em J. Phys. A}, 32(26):4867--4876, 1999.

\bibitem{jensen-guttmann-animals}
I.~Jensen and A.~J. Guttmann.
\newblock Statistics of lattice animals (polyominoes) and polygons.
\newblock {\em J. Phys. A}, 33(29):L257--L263, 2000.

\bibitem{klarner-results}
D.~A. Klarner.
\newblock Some results concerning polyominoes.
\newblock {\em Fibonacci Quart.}, 3:9--20, 1965.

\bibitem{klarner}
D.~A. Klarner.
\newblock Cell growth problems.
\newblock {\em Canad. J. Math.}, 19:851--863, 1967.

\bibitem{klarner-rivest}
D.~A. Klarner and R.~L. Rivest.
\newblock Asymptotic bounds for the number of convex {$n$}-ominoes.
\newblock {\em Discrete Math.}, 8:31--40, 1974.

\bibitem{borgne-marckert}
Y.~Le~Borgne and J.-F. Marckert.
\newblock Directed animals and gas models revisited.
\newblock {\em Electron. J. Combin.}, R71, 2007.

\bibitem{madras-slade}
N.~Madras and G.~Slade.
\newblock {\em The self-avoiding walk}.
\newblock Probability and its Applications. Birkh\"auser Boston Inc., Boston,
  MA, 1993.

\bibitem{polya}
G.~P{\'o}lya.
\newblock On the number of certain lattice polygons.
\newblock {\em J. Combinatorial Theory}, 6:102--105, 1969.

\bibitem{prellberg-brak}
T.~Prellberg and R.~Brak.
\newblock Critical exponents from nonlinear functional equations for partially
  directed cluster models.
\newblock {\em J. Stat. Phys.}, 78(3/4):701--730, 1995.

\bibitem{privman-barma}
V.~Privman and M.~Barma.
\newblock Radii of gyration of fully and partially directed animals.
\newblock {\em Z. Phys. B: Cond. Mat.}, 57:59--63, 1984.

\bibitem{svrakic}
V.~Privman and N.~M. {\v{S}}vraki{\'c}.
\newblock Exact generating function for fully directed compact lattice animals.
\newblock {\em Phys. Rev. Lett.}, 60(12):1107--1109, 1988.

\bibitem{prodinger}
H.~Prodinger.
\newblock The kernel method: a collection of examples.
\newblock {\em S\'em. Lothar. Combin.}, 50:Art. B50f, 19 pp. (electronic),
  2003/04.

\bibitem{read}
R.~C. Read.
\newblock Contributions to the cell growth problem.
\newblock {\em Canad. J. Math.}, 14:1--20, 1962.

\bibitem{rechni-SAP}
A.~Rechnitzer.
\newblock Haruspicy 2: the anisotropic generating function of self-avoiding
  polygons is not {D}-finite.
\newblock {\em J. Combin. Theory Ser. A}, 113(3):520--546, 2006.

\bibitem{richard-synthese}
C.~Richard.
\newblock Limit distributions and scaling functions.
\newblock Arxiv:0704.0716v2 [math-ph], 2007.

\bibitem{salomaa}
A.~Salomaa and M.~Soittola.
\newblock {\em Automata-theoretic aspects of formal power series}.
\newblock Springer-Verlag, New York, 1978.
\newblock Texts and Monographs in Computer Science.

\bibitem{temperley}
H.~N.~V. Temperley.
\newblock Combinatorial problems suggested by the statistical mechanics of
  domains and of rubber-like molecules.
\newblock {\em Phys. Rev. (2)}, 103:1--16, 1956.

\bibitem{viennot-heaps}
G.~X. Viennot.
\newblock Heaps of pieces. {I}. {B}asic definitions and combinatorial lemmas.
\newblock In {\em Combinatoire \'enum\'erative (Montr\'eal, 1985)}, volume 1234
  of {\em Lecture Notes in Math.}, pages 321--350. Springer, Berlin, 1986.

\bibitem{voege-guttmann}
M.~V\"oge and A.~J. Guttmann.
\newblock On the number of hexagonal polyominoes.
\newblock {\em Theoret. Comput. Sci}, 307(2):433--453, 2003.

\bibitem{yuba}
T.~Yuba and M.~Hoshi.
\newblock Binary search networks: a new method for key searching.
\newblock {\em Inform. Process. Lett.}, 24:59--65, 1987.

\bibitem{zeilberger-skinny}
D.~Zeilberger.
\newblock Symbol-crunching with the transfer-matrix method in order to count
  skinny physical creatures.
\newblock {\em Integers}, pages A9, 34pp. (electronic), 2000.

\bibitem{zeilberger-animals}
D.~Zeilberger.
\newblock The umbral transfer-matrix method. {III}. {C}ounting animals.
\newblock {\em New York J. Math.}, 7:223--231 (electronic), 2001.

\end{thebibliography}
 \end{document}